\magnification=\magstep1
\baselineskip 16pt

\centerline{\bf Commutative and Noncommutative Invariant Theory}
\bigskip
\centerline{\bf Vesselin Stoyanov Drensky
\footnote{$^{\ast,}$}{\rm
Partially supported by Grant MM404/94 of the Ministry of
Science and Education, Bulgaria.}
\footnote{$^{\ast\ast}$}{\rm Published as V. S. Drensky, Commutative and noncommutative invariant theory,
Math. and Education in Math., Proc. of the 24-th Spring Conf. of the Union of Bulgar. Mathematicians, Svishtov, April 4-7, 1995,  14-50.}}

\medskip

\centerline{\bf Abstract}
\medskip

The purpose of this survey paper is to bring to a large
mathematical audience (containing also non-algebraists) some topics
of invariant theory both in the classical commutative and the
recent noncommutative case. We have included only several topics
from the classical invariant theory -- the finite generating (the
Endlichkeitssatz) and the finite presenting (the Basissatz) of the
algebra of invariants, the Molien formula for its Hilbert series
and the Shephard-Todd-Chevalley theorem for the invariants of a
finite group generated by pseudo-reflections. Then we give
analogues of these results for free and relatively free associative
and Lie algebras. Finally we deal with the algebra of generic
matrices and the invariant theory related with it.
\medskip

\centerline{\bf Introduction}
\medskip

This is an extended version of the invited talk to be given at the
Annual Meeting of the Union of the Bulgarian Mathematicians. Our
purpose is to bring to a large mathematical audience (containing
also non-algebraists) some topics of invariant theory both in the
classical commutative and the recent noncommutative case. We do not
pretend for a comprehensive survey on the subject. We have chosen
several topics only. For a background on classical invariant theory
we refer to the books by Springer [87] or Dieudonn\'e and Carrell
[17]. We want also to mention the excellent survey by Formanek [36]
on noncommutative invariant theory which involved us in this topic
and inspired the present paper.

The paper is organized as follows. In Section 1 we give some
classical results on commutative invariant theory which are
considered to be corner stones of the theory. In Section 2 we
discuss analogues of these results when a finite linear group acts
on the free associative algebra and in Section 3 we extend our
considerations to relatively free algebras. Section 4 is devoted to
the case of finite linear groups acting on relatively free Lie
algebras. Finally Section 5 deals with invariants of the projective
special linear group $PSL_k$ acting simultaneously by conjugation
on several $k\times k$ matrices.

\medskip
\centerline{\bf 1. Classical Invariant Theory}
\medskip

Let $G$ be a group acting on a set $V$, i.e. there is a map
$$
G \times V \longrightarrow V,\, \hbox{\rm which we denote by}\,
(g,v) \longrightarrow g(v),
$$
such that

\item{(i)} $(gh)(v) = g(h(v))$ for any $g,h \in G$ and $v \in V$;
\par
\item{(ii)} $1(v) = v$ for any $v \in V$ and 1 is the unit element
of $G$.
\par
\noindent The element $v \in V$ is called an {\bf invariant} of $G$
if $g(v) = v$ for any $g\in G$. We denote the set of all invariants
of $G$ by $V^G$.

If the group $G$ acts on the sets $V$ and $W$ we associate an
action of $G$ on the set ${\cal F}(V,W)$ of all functions from $V$
to $W$ by
$$
(g(f))(v) = f(g^{-1}(v)),\, f \in {\cal F}(V,W),\, g \in G, v \in V.
$$
The function $f$ is called a {\bf concominant} of $G$ if it is
invariant for the action of $G$ on ${\cal F}(V,W)$. Equivalently
$f\in {\cal F}(V,W)$ is a concominant if $g(f(v)) = f(g(v))$ for
all $g \in G$, $v \in V$.

{\sl All the considerations in this paper are over a fixed
algebraically closed field $K$ of characteristic $0$.} Let $V_m$ be
a vector space with basis $v_1,\ldots,v_m$. We assume that the
general linear group $GL_m = GL(V_m) = GL_m(K)$ acts canonically on
$V_m$ and identify it with the group of $m\times m$ invertible
matrices. If the group $G$ acts on $V_m$ it is naturally to require
for this action to be linear. It induces a homomorphism $\rho:G
\longrightarrow GL_m$; $\rho$ is called a {\bf representation} of
$G$ in $V_m$ and $V_m$ is a $G$-{\bf module}. The function
$\chi_{\rho}:G\longrightarrow K$ defined by $\chi_{\rho}(g) =
{\rm tr}_{V_m}(\rho(g))$, $g\in G$, is called the {\bf character}
of $\rho$. The representation $\rho$ is {\bf completely reducible}
if $V_m$ is a direct sum of minimal $G$-invariant subspaces. If the
group $G$ is finite the {\bf Maschke theorem} states that every
representation of $G$ in $V_m$ is completely reducible. For our
purposses we may assume without loss of generality that $G$ is a
subgroup of $GL_m$ and consider the elements of $G$ as
invertible $m\times m$ matrices. Let $\psi:G\longrightarrow GL_p$
be another representation of $G$. If the entries $\psi_{pq}(g)$ of
the matrix $\psi(g)$ are rational (respectively polynomial)
functions of the entries $g_{ij}$ of the matrices $g\in G$ the
representation $\psi$ is called {\bf rational} (respectively {\bf
polynomial}).

Let $X_m =\{x_1,\ldots,x_m\}$ and let ${\rm span}(X_m)$ be the
vector space with basis $x_1,\ldots,x_m$. The polynomial algebra
$K[X_m] = K[x_1,\ldots,x_m]$ in $m$ commuting variables
$x_1,\ldots,x_m$ can be viewed as the algebra of polynomial
functions on $V_m$. Its generators $x_i$ are defined via
$$
x_i(\xi_1v_1 + \ldots + \xi_mv_m) = \xi_i,\, \xi_j \in K,\, i,j =
1,\ldots,m.
$$
Equipping $K$ with the trivial action of $G$ the above arguments
give that the action of $GL_m$ on $V_m$ induces an action on
$K[x_1,\ldots,x_m]$ by
$$
(gf(x_1,\ldots,x_m))(v) = f(x_1,\ldots,x_m)(g^{-1}(v)),
$$
where $g\in GL_m$, $f(x_1,\ldots,x_m)\in K[x_1,\ldots,x_m]$ and
$v\in V_m$. In this way ${\rm span}(X_m)$ is isomorphic to the dual
space of $V_m$ and $GL_m$ is embedded into the group of
automorphisms of the polynomial algebra. Replacing $V_m$ with its
dual space {\sl till the end of this paper we shall assume that
$GL_m$ acts canonically on ${\rm span}(X_m)$} (and not on
$V_m$). {\sl If $g =(g_{ij})\in GL_m$ then $g(x_j) = \sum_{i=1}^m
g_{ij}x_i$, $j = 1,\ldots,m$.}

{\bf Definition 1.1.} {\sl Let $G$ be a subgroup of $GL_m$. The
function $f(x_1,\ldots,x_m) \in K[x_1,\ldots,x_m]$ is called a
({\bf polynomial}) {\bf invariant} of $G$ if $f$ is a concominant
of $G$, i.e.
$$
gf(x_1,\ldots,x_m) = f(x_1,\ldots,x_m)\, \hbox{\rm for every}\, g
\in G.
$$
The {\bf algebra of invariants} of $G$ is the subalgebra of
$K[x_1,\ldots,x_m]$}
$$
K[x_1,\ldots,x_m]^G = \{f \in K[x_1,\ldots,x_m] \mid f\, \hbox{\rm
is invariant of}\, G\}.
$$

The origins of the invariant theory can be found in the 18th
century. One of its main problems is to describe for a fixed group
$G$ the algebra of invariants of $G$. In the 19th century the first
methods were developed to attack this problem. The modern approach
begins with the results of Hilbert [47, 48] and Emmy Noether [71].
In particular Hilbert and Noether have established the following
famous theorems.

{\bf Theorem 1.2.} (Endlichkeitssatz [71]) {\sl Let $G$ be a
finite subgroup of $GL_m$. Then the algebra of invariants
$K[x_1,\ldots,x_m]^G$ is finitely generated. It has a system of
generators $f_1,\ldots,f_p$ where every $f_i$ is a homogeneous
polynomial of degree bounded by the order $\vert G\vert$ of the
group $G$.}

{\bf Theorem 1.3.} (Basissatz [47]) {\sl Every commutative algebra
$S$ with a finite set of generators $f_1,\ldots,f_p$ can be defined
by a finite system of relations, i.e. there exists a finite system
of polynomials
$$
r_1(y_1,\ldots,y_p),\ldots,r_q(y_1,\ldots,y_p)\in K[y_1,\ldots,y_p]
$$
such that the kernel of the canonical homomorphism
$$
\nu: K[y_1,\ldots,y_p] \longrightarrow S,\,\hbox{\sl where}\,
\nu(y_j) = f_j,\, j = 1,\ldots,p,
$$
is generated as an ideal by
$r_1(y_1,\ldots,y_p),\ldots,r_q(y_1,\ldots,y_p)$.}

Theorem 1.2 can be also generalized for some classes of infinite
groups. Here we give a simplified version of a result of Nagata
which extends a theorem of Hilbert.

{\bf Theorem 1.4.} (Hilbert-Nagata, see e.g. [17]) {\sl Let $G$
be a subgroup of $GL_m$ such that any finite dimensional rational
representation of $G$ is completely reducible. Then the algebra of
invariants $K[x_1,\ldots,x_m]^G$ is finitely generated.}

Between the groups covered by Theorem 1.4 are the finite groups,
$GL_p$, the {\bf special linear group} $SL_p = \{ g\in GL_p\mid
{\rm det}(g) = 1\}$, all semisimple algebraic groups. In his
famous lecture {\sl ``Mathematische Probleme''} given at the
International Congress of Mathematicians held in 1900 in Paris
Hilbert [49] asked as Problem 14 {\sl if the algebra
$K[x_1,\ldots,x_m]^G$ is finitely generated for every subgroup $G$
of $GL_m$.} A negative answer was given by Nagata [70] (see also
[17]).

Traditionally a result giving the explicit generators of the
algebra of invariants of a group $G\subset GL_m$ is called the {\bf
First Fundamental Theorem of the Invariant Theory} of $G$ and a
result describing the relations between the generators -- the {\bf
Second Fundamental Theorem}.

For example, if the symmetric group $S_m$ acts canonically on the
basis of ${\rm span}(X_m)$ by $\sigma(x_i) = x_{\sigma(i)}$,
$\sigma\in S_m$, $i = 1,\ldots,m$, the First Fundamental Theorem
states that the elementary symmetric functions
$e_i(x_1,\ldots,x_m)$, $i =
1,\ldots,m$, generate the algebra of symmetric polynomials in $m$
variables and the Second Fundamental Theorem says that there are no
relations between them, i.e. they are algebraically independent.
The complete description of the finite groups $G$ such that the
algebra of invariants $K[x_1,\ldots,x_m]^G$ is isomorphic to a
polynomial algebra is given by Shephard and Todd [84] and Chevalley
[15] (see also [87]). Recall that a matrix $g\in GL_m$ is called a
{\bf pseudo-reflection} if $n - 1$ of its eigenvalues are equal to 1
and if moreover ${\rm span}(X_m)$ has a basis of eigenvectors of
$g$.

{\bf Theorem 1.5.} (Shephard-Todd-Chevalley [84, 15])
{\sl For a finite subgroup $G$ of
$GL_m$ the algebra of invariants $K[x_1,\ldots,x_m]^G$ is
isomorphic to a polynomial algebra if and only if $G$ is generated
by pseudo-reflections.}

The vector space $W$ is {\bf graded} if it is a direct sum of its
subspaces $W^{(n)}$, $n = 0,1,2,\ldots\,$; the subspace $W^{(n)}$
is the {\bf homogeneous component of degree} $n$ of $W$. In our
considerations all $W^{(n)}$ are finite dimensional. The formal
power series
$$
{\rm Hilb}(W,t) = \sum_{n=0}^{\infty}({\rm dim}_KW^{(n)})t^n
$$
is called the {\bf Hilbert} (or {\bf Poincar\'e}) {\bf series} of
$W$. Traditionally if the Hilbert series is convergent we identify
it with another function $f(t)$ if ${\rm Hilb}(W,t) = f(t)$ in a
neighbourhood of 0. Similarly the vector space $W$ is {\bf
multigraded} if it is a direct sum of its {\bf multihomogeneous
components} $W^{(n_1,\ldots,n_m)}$ and the {\bf Hilbert series} of
$W$ is the formal power series in $m$ variables
$$
{\rm Hilb}(W,t_1,\ldots,t_m) =
\sum_{n_1=0}^{\infty}\ldots\sum_{n_m=0}^{\infty}{\rm
dim}_KW^{(n_1,\ldots,n_m)}t_1^{n_1}\ldots t_m^{n_m}.
$$
For example the polynomial algebra $K[x_1,\ldots,x_m]$ is
canonically graded assuming that the monomials of degree $n$ span
its homogeneous component of degree $n$. Similarly we may define a
multigrading on $K[x_1,\ldots,x_m]$ counting the degree of each
variable. Clearly the multihomogeneous components of
$K[x_1,\ldots,x_m]$ are one-dimensional and
$$
{\rm Hilb}(K[x_1,\ldots,x_m],t_1,\ldots,t_m) =
\sum_{n_1=0}^{\infty}\ldots\sum_{n_m=0}^{\infty}t_1^{n_1}\ldots
t_m^{n_m} =
\prod_{i=1}^m{1\over 1 - t_1}.
$$

For any group $G\subseteq GL_m$ the algebra of invariants of $G$
inherits the grading of $K[x_1,\ldots,x_m]$ and the Molien formula
gives its Hilbert series for finite groups.

{\bf Theorem 1.6.} (Molien formula [68])
{\sl The Hilbert series of the algebra of
invariants of a finite subgroup $G$ of $GL_m$ is}
$$
{\rm Hilb}(K[x_1,\ldots,x_m]^G,t) =
{1\over \vert G\vert}\sum_{g\in G}{1\over {\rm det}(1 - gt)}.
$$
\medskip

\centerline{\bf 2. Free Associative Algebras}
\medskip

As in the previous section ${\rm span}(X_m)$, $X_m
=\{x_1,\ldots,x_m\}$, is the $m$-dimensional vector space with
basis $x_1,\ldots,x_m$ with the canonical action of $GL_m$. The
{\bf free associative algebra} $K\langle X_m\rangle = K\langle
x_1,\ldots,x_m\rangle$ {\bf of rank} $m$ freely generated by the
set $X_m = \{x_1,\ldots,x_m\}$ (or the algebra of polynomials in
$m$ noncommuting variables) is the vector space with basis all the
monomials $x_{i_1}\ldots x_{i_n}$, $i_k = 1,\ldots,m$, $k =
1,\ldots,n$, $n = 0,1,2,\ldots,$ and multiplication defined on the
basis by
$$
(x_{i_1}\ldots x_{i_p})(x_{j_1}\ldots x_{j_q}) =
x_{i_1}\ldots x_{i_p}x_{j_1}\ldots x_{j_q}.
$$
The algebra $K\langle X_m\rangle$ is also defined by the
following universal property which is similar to the universal
property of the ordinary polynomial algebra in the class of all
commutative algebras: $K\langle X_m\rangle$ is generated by the
set $X_m$ and for any associative algebra $R$ every map
$X_m\longrightarrow R$ can be uniquely extended to a homomorphism
$K\langle X_m\rangle\longrightarrow R$. It is easy to see that
$K\langle X_m\rangle$ is isomorphic to the tensor algebra of the
vector space ${\rm span}(X_m)$. The action of $GL_m$ on ${\rm
span}(X_m)$ is extended diagonally on $K\langle X_m\rangle$, i.e.
$$
g(x_{i_1}\ldots x_{i_n}) = g(x_{i_1})\ldots g(x_{i_n}),\,g\in
GL_m.
$$
The free algebra $K\langle X_m\rangle$ has a natural multigrading
which is similar to that of the polynomial algebra $K[X_m]$ and
the monomials of degree $n_i$ with respect to $x_i$ span the
multihomogeneous component of $K\langle X_m\rangle$ of multidegree
$(n_1,\ldots,n_m)$; the canonical grading of $K\langle X_m\rangle$
is defined counting the total degree of the monomials. Clearly the
homogeneous component $K\langle X_m\rangle^{(n)}$ has a basis
consisting of all monomials $x_{i_1}\ldots x_{i_n}$, $i_j =
1,\ldots,m$, $j = 1,\ldots,n$,
and ${\rm dim} K\langle X_m\rangle^{(n)} = m^n$.
Hence
$$
{\rm Hilb}(K\langle X_m\rangle,t) = \sum_{n=0}^{\infty} m^nt^n =
{1\over 1 - mt}.
$$
Similarly we obtain that
$$
{\rm Hilb}(K\langle X_m\rangle^{(n_1,\ldots,n_m)},t_1,\ldots,t_m)
=\sum_{i_1=1}^m\ldots\sum_{i_n=1}^mt_{i_1}\ldots t_{i_n} = (t_1 +
\ldots + t_m)^n
$$
and
$$
{\rm Hilb}(K\langle X_m\rangle^{(n_1,\ldots,n_m)},t_1,\ldots,t_m)
={1\over 1 -(t_1 +\ldots +t_m)}.
$$
In this section we discuss analogues of the results of Section 1
to the case when a finite linear group acts on $K\langle
X_m\rangle$.

{\bf Definition 2.1.} {\sl The group $G \subset GL_m$ acts on
${\rm span}(X_m)$ by {\bf scalar multiplication} if $G$ consists of scalar
matrices.}

Clearly if $G$ is finite and acts by scalar multiplication, then
$G$ is cyclic; if $\vert G\vert = d$ then $K\langle X_m\rangle^G$
is isomorphic to the tensor algebra of the vector space $K\langle
X_m\rangle^{(d)}$ of all homogeneous noncommutative polynomials of
degree $d$. Since ${\rm dim} K\langle X_m\rangle^{(d)} = m^d$ the
algebra $K\langle X_m\rangle^G$ is isomorphic to the free
associative algebra of rank $m^d$. It turns out that the analogue
of the Emmy Noether Theorem 1.2 holds for $K\langle X_m\rangle^G$
in this very special case only.

{\bf Theorem 2.2.} (Korjukin, Dicks and Formanek [16], Kharchenko
[56]) {\sl Let $G$ be a finite subgroup of $GL_m$. Then $K\langle
X_m\rangle^G$ is finitely generated if and only if $G$ acts on
${\rm span}(X_m)$ by scalar multiplication.}

The analogue of the Shephard-Todd-Chevalley Theorem 1.5 sounds
also very different for $K\langle X_m\rangle$. It turns out that
the algebra $K\langle X_m\rangle^G$ is always free. Additionally
there is a Galois correspondence between the subgroups of $G$ and
the free subalgebras of $K\langle X_m\rangle$ which contain
$K\langle X_m\rangle^G$.

{\bf Theorem 2.3.} (Lane [60], Kharchenko [55]) {\sl For every
finite subgroup G of GL  the algebra of invariants $K\langle
X\rangle^G$ is free.}

{\bf Theorem 2.4.} (Kharchenko [55]) {\sl The map
$H\longrightarrow K\langle X_m\rangle^H$ gives a 1-1
correspondence between the subgroups of $G$ and the free
subalgebras of $K\langle X_m\rangle$ containing $K\langle
X_m\rangle^G$.}

Finally the Molien formula from Theorem 1.6 has the following nice
form.

{\bf Theorem 2.5.} (Dicks and Formanek [16]) {\sl The Hilbert
series of the algebra of invariants $K\langle X_m\rangle^G$ of a
finite subgroup $G$ of $GL_m$ is}
$$
{\rm Hilb}
(K\langle X_m\rangle^G,t) = {1\over\vert G\vert}\sum_{g\in
G}{1\over1 - {\rm tr}(g)t}.
$$

\medskip
\centerline{\bf 3. Groups Acting on Relatively Free Algebras}
\medskip

If $I$ is an ideal of $K\langle X_m\rangle$ which is
$GL_m$-invariant, i.e. $GL_m(I) = I$, then every $g \in GL_m$
induces an automorphism of the factor-algebra $K\langle
X_m\rangle/I$ and we may consider the algebra of invariants
$(K\langle X_m\rangle/I)^G$ for any subgroup $G$ of $GL_m$.
Especially interesting are the ideals $I$ which are stable under
all the endomorphisms of $K\langle X_m\rangle/I$. Such ideals have
another description related with the theory of PI-algebras.

{\bf Definition 3.1.} {\sl Let $R$ be an associative algebra.}

(i) {\sl A polynomial $f(x_1,\ldots,x_m) \in K\langle
x ,\ldots,x_m\rangle$ is called a {\bf polynomial identity} for
$R$ if $f(r_1,\ldots,r_m) = 0$ for any $r_1,\ldots,r_m \in R$. If
$R$ satisfies a non-trivial polynomial identity $R$ is called a
{\bf PI-algebra}. The class of all algebras satisfying a given
system of polynomial identities $\{f_j(x_1,\ldots,x_{n_j})\mid j
\in J\}$ is the {\bf variety} of algebras determined by this
system. If the variety is determined by all the polynomial
identities of $R$ it is called the variety {\bf generated} by $R$
and is denoted by ${\rm var}(R)$.}

(ii)  {\sl The set $T_m(R)$ of all polynomial identities in $m$
variables for $R$ is called a {\bf T-ideal} of $K\langle
x_1,\ldots,x_m\rangle$. The factor-algebra $F_m({\rm var} R) =
F_m(R) = K\langle X_m\rangle/T_m(R)$ is the {\bf relatively free
algebra of rank} $m$ in the variety ${\rm var}(R)$.}

Typical examples of PI-algebras are the commutative algebras
satisfying the {\bf commutator identity}
$$
[x_1,x_2] = x_1x_2 - x_2x_1 = 0
$$
and the finite dimensional algebras. If ${\rm dim} R < n$ then $R$
satisfies the {\bf standard identity} of degree $n$
$$
s_n(x_1,\ldots,x_n) = \sum_{\pi\in S_n}({\rm sign}
\pi)x_{\pi(1)}\ldots x_{\pi(n)} = 0.
$$
A non-trivial example of a PI-algebra is the {\bf Grassmann} (or
{\bf exterior}) {\bf algebra} $E(V)$ of a vector space $V$. If $V$
has a basis $v_1,v_2,\ldots$, then $E(V)$ has a basis
$v_{i_1}\ldots v_{i_n}$, where $i_1 < \ldots < i_n$, $n =
0,1,2,\ldots$ and defining relations $e_ie_j = -e_je_i$, $i,j =
1,2,\ldots\, .$ It is easy to see that for any elements $u_1,u_2
\in E(V)$, the commutator $[u_1,u_2]$ belongs to the subspace of
$E(V)$ spanned on the words $v_{i_1} \ldots v_{i_{2n}}$ of even
length and is in the centre of $E(V)$. Therefore $E(V)$ satisfies
the polynomial identity
$$
[x_1,x_2,x_3] = [[x_1,x_2],x_3] = 0.
$$
Krakowski and Regev [58] have shown that if ${\rm dim} V = \infty$
then all the polynomial identities of $E(V)$ follow from
$[x_1,x_2,x_3] = 0$, i.e. $T_{\infty}(E(V)) \subset K\langle
x_1,x_2,\ldots\rangle$ coincides with the minimal T-ideal
containing $[x_1,x_2,x_3]$.

It is easy to see that an ideal $I$ of $K\langle X_m\rangle$ is
stable under all the endomorphisms of $K\langle X_m\rangle$ if and
only if $I = T_m(R)$ for some PI-algebra $R$. An example of such an
algebra $R$ with $T_m(R) = I$ is $R = K\langle X_m\rangle/I$. Old
and recent results on structure ring theory show that the class of
PI-algebras enjoys a lot of good properties of finite dimensional
and commutative algebras (see e.g. the books by Jacobson [50],
Rowen [82] and Kemer [54]). The investigation of the relatively
free algebras is
stimulated not only by invariant theory. All the information about
the polynomial identities for a PI-algebra R may be obtained from
its relatively free algebras $F_m(R)$, $m = 1,2,\ldots\, .$If $R$
is finite dimensional it is sufficient to study $F_m(R)$ for $m =
1,2,\ldots,{\rm dim} R$ only.

Every T-ideal $T_m(R)$ of $K\langle X_m\rangle$ is a
multihomogeneous vector subspace of $K\langle X_m\rangle$, i.e.
$T_m(R) = \sum^{\oplus}((K\langle X_m\rangle)^{(n_1,\ldots,n_m)}
\cap T_m(R))$, and the relatively free algebra $F_m(R)$ inherits
the (multi)grading of $K\langle X_m\rangle$. Similarly $F_m(R)$
inherits the action of $GL_m$ on $K\langle X_m\rangle$. It turns
out that the Hilbert series ${\rm Hilb}(F_m(R),t_1,\ldots,t_m)$
plays the r\^ole of the character of the representation of $GL_m$ in
$F_m(R)$.

{\bf Theorem 3.2.} (Formanek [35, Sections 2 and 3]) {\sl Let $F =
F_m(R)$ be a relatively free algebra with Hilbert series ${\rm
Hilb}(F_m(R),t_1,\ldots,t_m)$ and let $g = \xi_1e_{11} + \ldots +
\xi_me_{mm} \in GL_m$ be a diagonal matrix with entries
$\xi_1,\ldots,\xi_m$ on the diagonal. Then the formal power series
${\rm Hilb}(F_m(R),\xi_1t,\ldots,\xi_mt)$ is equal to the
generating function of the sequence ${\rm tr}_{F^{(n)}}(g)$,
$n = 0,1,2,\ldots$, i.e.}
$$
\chi_{F^{(n)}}(g) = {\rm tr}_{F^{(n)}}(g) =
\sum_{n_1+\ldots+n_m=n}{\rm dim}
F^{(n_1,\ldots,n_m)}\xi_1^{n_1}\ldots\xi_m^{n_m}.
$$

In all the cases of PI-algebras when the explicit form of the
Hilbert series of $F_m(R)$ is known the computing of ${\rm
Hilb}(F_m(R),t_1,\ldots,t_m)$ is a consequence of the description
of the $GL_m$-module structure of $F_m(R)$. Traditionally (see the
results of Regev, e.g. [81]) one works with the sequence
$\chi_{S_n}(R)$, $n = 0,1,2,\ldots$, of the $S_n$-{\bf
cocharacters} of an algebra $R$. It is defined as the sequence of
characters of the symmetric group $S_n$ acting on the vector space
$P_n/(P_n \cap T_n(R))$, $n = 0,1,2,\ldots$, where $P_n$ is the set
of all {\bf multilinear} in $x_1,\ldots,x_n$ {\bf polynomials} from
$K\langle X_n\rangle$. It is well known that the irreducible
representations of $S_n$ and the irreducible polynomial
representations of $GL_m$ are described in the same way by
partitions and Young diagrams. The author [19] and Berele [7] have
independently shown that the $S_n$-module structure of $P_n/(P_n
\cap T_n(R))$ and the $GL_m$-module structure of $F_m^{(n)}(R)$ are
the same for any $n = 0,1,2,\ldots\, .$

Translating the results of [58] to the language of $GL_m$-modules,
in the special case of the Grassmann algebra $E$ of an infinite
dimensional vector space we obtain (see e.g. [14])
$$
{\rm Hilb}(F_m(E),t_1,\ldots,t_m) =\prod_{i=1}^m{1+t_i\over 1-t_i}.
$$
From the point of view of structure theory of T-ideals the most
important PI-algebras are the matrix algebras $M_k(K)$, the matrix
algebras $M_k(E)$ with entries from the Grassmann algebra $E$ and
some subalgebras of $M_k(E)$ (see [54]). In particular by a
classical result of Amitsur (see [82, Corollary 2.4.10, p. 135])
$F_m(M_k(K))$ are the only semiprime relatively free algebras. As a
consequence of the theorem of Razmyslov-Kemer-Braun for the
nilpotency of the radical of finitely generated PI-algebras [79,
52, 11] it follows that for any PI-algebra $R$ and any positive
integer $m$ there exists a pair of integers $k$ and $c$ such that
$T_m(M_k(K)) \supseteq T_m(R) \supseteq T_m^c(M_k(K))$. In this way
the polynomial identities of the matrix algebras may be used as a
measure how complicated are the polynomial identities in $m$
variables for an algebra $R$. More precise classification counting
all the polynomial identities for $R$ has been given by Kemer [54] who
has developed the structure theory of T-ideals in the spirit of the
structure theory of ideals in the polynomial algebra. From this
point of view the simplest T-ideals are the non-matrix ones, i.e.
these which are not contained in $T_m(M_2(K))$. By a result of
Latyshev [61] if $R$ is a finitely generated PI-algebra with a
non-matrix T-ideal then $R$ satisfies the polynomial identity
$[x_1,x_2] \ldots [x_{2c-1},x_{2c}] = 0$ for some $c$. In Section 5
we shall give another realization of $F_m(M_k(K))$ as the algebra
of generic $k \times k$ matrices. The explicit expression for the
Hilbert series of $F_m(M_k(K))$ is known in the case of $2 \times
2$ matrices only (see [31] for $m = 2$ and [35, 76, 20] for
$m$ arbitrary).

Now we survey the analogues for relatively free algebras of the
main results of Section 1. A priori it is difficult to guess the
behaviour of the algebra of invariants in this case. Many results
on combinatorial and structure ring theory show that the relatively
free algebras are closer to the polynomial algebra than to the free
associative algebra. Nevertheless if the T-ideal $T_m(R)$ is ``big
enough'' it is naturally to expect that many properties of
$K[X_m]^G$ can be transferred to $F_m(R)^G$ and if $T_m(R)$ is
``small'' the behaviour of $F_m(R)^G$ is as of $K\langle
X_m\rangle^G$. We start with the problem of finite generating of
$F_m(R)^G$. For different PI-algebras $R$ and different finite
groups $G$ the answer should varies somewhere between Theorem 1.2
and Theorem 2.2. It turns out that the T-ideals $T_m(R)$ such that
the algebra $F_m(R)^G$ is finitely generated for any group $G$ have
several different descriptions in the languages of noetherian
algebras, polynomial identities (also for $2 \times 2$ matrices)
and even decision problems at the meeting of algebra and logic.

First we define a sequence of PI-algebras $R_k$,
$k = 1,2,\ldots\, .$ Let $C_k = K[t]/(t^k)$ and let
$$
R_k = \pmatrix{C_k&tC_k\cr
tC_k&C_k\cr}
\subset M_2(C_k),
$$
where $M_2(C_k)$ is the $2 \times 2$ matrix algebra with entries
from $C_k$. These algebras ``almost'' describe
the T-ideals containing strictly the T-ideal $T_m(M_2(K))$ [21]. If
$T_m(R) \supset T_m(M_2(K))$, $m > 1$, then there exist a positive
integer $k$ and a finite dimensional algebra $N$ such that $N$ is
one-dimensional modulo its Jacobson radical and $T_m(R) = T_m(R_k
\oplus N)$. Another description of these T-ideals is given in
[53].

{\bf Theorem 3.3.} {\sl Let $R$ be a PI-algebra. The following
conditions on the algebra $R$ are equivalent. If some of them is
satisfied for some $m_0 \geq 2$, then all of them hold for all $m
\geq 2$:}

(i) {\sl The algebra $F_m(R)^G$ is finitely generated for every
finite subgroup $G$ of $GL_m$.}

(ii) {\sl The algebra $F_m(R)^{\langle g\rangle}$ is finitely
generated, where $g \in GL_m$ is a matrix of finite multiplicative
order with at least two eigenvalues (or characteristic roots) of
different order.}

(iii) {\sl The algebra $F_m(R)$ is {\bf weakly noetherian}, i.e.
satisfies the ascending chain condition for two-sided ideals.}

(iv) {\sl Let $S$ be an algebra satisfying all the polynomial
identities of $R$ (i.e. $S \in {\rm var}(R))$ and generated by $m$
elements $s_1,\ldots,s_m$. Then $S$ is {\bf finitely presented} as
a homomorphic image of $F_m(R)$, i.e. the kernel of the canonical
homomorphism $F_m(R) \longrightarrow S$ defined by $x_i
\longrightarrow s_i$, $i = 1,\ldots,m$, is a finitely generated
ideal of $F_m(R)$.}

     (v) {\sl If $S$ is a finitely generated algebra from ${\rm
var}(R)$, then $S$ is {\bf residually finite}, i.e. for every
element $s \in S$ there exist a finite dimensional algebra D and a
homomorphism $\phi:S \longrightarrow D$ such that $\phi(s) \not=
0$.}

(vi) {\sl If $S$ is a finitely generated algebra from ${\rm
var}(R)$, then $S$ is {\bf representable by matrices}, i.e. there
exist an extension $L$ of the base field $K$ and an integer $k$
such that $S$ is
isomorphic to a subalgebra of the $K$-algebra $M_k(L)$ of all
$k\times k$ matrices with entries from $L$.}

(vii) {\sl If the base field $K$ is countable, then  every set of
pairwise non-isomorphic homomorphic images of $F_m(R)$ is
countable.}

(viii) {\sl The algebra $R$ satisfies a polynomial identity of the
form}
$$
x_2x_1^{n-2}x_2 +
\sum_{i+j>0}\alpha_{ij}x_1^ix_2x_1^{n-i-j-2}x_2x_1^j = 0,\,
\alpha_{ij} \in K.
$$

(ix) {\sl The algebra $R$ satisfies the polynomial identity
$$
[x_1,x_2,\ldots,x_2]x_3^n[x_4,x_5,\ldots,x_5] = 0
$$
for sufficiently long commutators and $n$ large enough.}

(x) {\sl If the algebra $R$ is finitely generated then it satisfies
the polynomial identity
$$
[x_1,\ldots,x_n]x_{n+1}\ldots x_{2n}[x_{2n+1},\ldots,x_{3n}] = 0
$$
for some positive integer $n$.}

(xi) {\sl The algebra $R$ satisfies a polynomial identity which
does not follow from the polynomial identities}
$$
[x_1,x_2][x_3,x_4][x_5,x_6] = 0,
$$
$$
[[x_1,x_2][x_3,x_4],x_5] = 0,
$$
$$
s_4(x_1,x_2,x_3,x_4) = 0.
$$

(xii) {\sl The T-ideal $T_m(R)$ is not contained in the T-ideal
$T_m(R_3)$ defined above.}

The equivalence of (i) and (iii) is established by Kharchenko [56],
of (iii) and (viii) by L'vov [65], of (iii), (v), (vi), (ix) and
(x) by Anan'in [3], of (viii) and (xi) by Tonov [92], of (ii),
(viii) and (xii) by the author [23] and the equivalence of (iii),
(iv) and (vii) is obvious. We should recall the paper by Malcev [66] which
was the starting point for many of the results on representable
algebras. The condition (ii) is a generalization of a result of
Fisher and Susan Montgomery [30] which states the following. If
$T_m(R) \subseteq T_m(M_2(K))$, $g \in GL_m$, $g^n = 1$ and $g$ has
at least two characteristic roots of different multiplicative
order, then the algebra of invariants $F_m(R)^{\langle g\rangle}$
is not finitely generated. It seems to us that (ii) may be used as
a simple criterion to check the equivalent conditions of Theorem
3.3. It is sufficient to choose $m = 2$ and
$$
g = \pmatrix{-1&0\cr
0&1\cr}.
$$
If $F_2(R)^{\langle g\rangle}$ is finitely generated then all the
assertions (i) - (xii) hold for $F_m(R)$ and $T_m(R)$ for any $m >
1$.

As we have mentioned above, the most attractive relatively free
algebras to study are $F_m(M_k(K))$. We give two results of
Formanek and Schofield [42] and Fisher and Susan Montgomery [30].

{\bf Theorem 3.4.} (i) (Formanek and Schofield [42]) {\sl Let $G$
be a finite subgroup of the group $SL_2$ of all $2 \times 2$
matrices with determinant $1$. Then the algebra of invariants
$F_2(M_2(K))^G$ is finitely generated.}

(ii) (Fisher and Susan Montgomery [30]) {\sl Let $G$ be a finite
cyclic subgroup of $GL_m$ which does not act on ${\rm span}(X_m)$
by scalar multiplication. Then there exists an integer $k$ such
that the algebra of invariants $F_m(M_k(K))^G$ is not finitely
generated.}

In virtue of Theorem 3.3 (ii) studying the algebra $F_m(R)^{\langle
g\rangle}$ it is important if the eigenvalues of $g \in GL_m$ are of
different multiplicative order. In the simplest case of algebras
with non-matrix T-ideals we have the following result.

{\bf Theorem 3.5.} [23] {\sl Let $R$ be a PI-algebra which
satisfies a non-matrix polynomial identity (i.e. $T_n(R)$ is not
contained in $T_n(M_2(K))$ for some $n$). If $g \in GL_m$ is a
matrix with all the eigenvalues of the same order then the algebra
of invariants $F_m(R)^{\langle g\rangle}$ is finitely generated.}

We are not going to discuss the problem of finite generation of
$F_m(R)^G$ (compare with the Hilbert Basissatz of Theorem 1.3). By
a result of Markov [67] very few T-ideals $T_m(R)$ are finitely
generated as ideals of $K\langle X_m\rangle$.

{\bf Theorem 3.6.} (Markov [67]) {\sl A T-ideal $T_m(R)$, $m > 1$,
is finitely generated as an ideal if and only if the algebra $R$
satisfies the Engel identity $x_2({\rm ad} x_1)^n = 0$ for some
$n$.}

Here $x_2({\rm ad} x_1) = [x_2,x_1]$ is the operator of the Lie
multiplication of $K\langle X_m\rangle$.

Very important is the result of Kemer [54] which gives the finite
generation of the T-ideals as T-ideals and answers into
affirmative the famous Specht problem [86].

{\bf Theorem 3.7.} (Kemer [54]) {\sl Every T-ideal $T(R)$ of the
free associative algebra $K\langle X_{\infty}\rangle = K\langle
x_1,x_2,\ldots\rangle$
of infinite rank is finitely generated as a T-ideal, i.e. the T-ideals
of $K\langle X_{\infty}\rangle$ (and hence in $K\langle
X_m\rangle$) satisfy the ascending chain condition.}

It turns out that the analogue of the Chevalley-Shephard-Todd
Theorem 1.5 holds very rarely for relatively free algebras. Let $G$
be a finite nontrivial subgroup of $GL_m$. As a consequence of the
inclusion $T_m(M_k(K)) \supseteq T_m(R) \supseteq T_m^c(M_k(K))$
for some $k$ and $c$ it follows that $T_m(M_k(K))/T_m(R)$ is the
Jacobson radical of $F_m(R) \cong K\langle X_m\rangle/T_m(R)$ and
the algebra of invariants $F_m(R)^G$ is relatively free if and only
if the algebra $F_m(M_k(K))^G$ is relatively free. Guralnick [43]
has established that $F_m(M_k(K))^G$ is not relatively free for any
$k > 1$ and $m > 1$. Therefore the algebras $R$ with $F_m(R)^G$
relatively free should satisfy a non-matrix polynomial identity.
The complete answer has been recently given by Domokos [18].

{\bf Theorem 3.8.} (Guralnick [43], Domokos [18]) {\sl For a finite
nontrivial subgroup $G$ of $GL_m$, $m > 1$, and a PI-algebra $R$,
the algebra of invariants $F_m(R)^G$ is relatively free if and only
if $G$ is generated by pseudo-reflections and $R$ satisfies the
polynomial identity $[x_2,x_1,x_1] = 0$.}

By [58] the polynomial identity $[x_1,x_2,x_3] = 0$ generates the
T-ideal of the Grassmann algebra $E$ of an infinite dimensional
vector space and it is easy to see that the polynomial identities
$[x_1,x_2,x_3] = 0$ and $[x_2,x_1,x_1] = 0$ are equivalent.
Therefore the analogue of Theorem 1.5 holds in two extremal cases
only -- either the algebra $R$ is very close to commutative or it
satisfies no polynomial identity (see Theorem 2.2).

Finally the following result of Formanek [36] is an analogue for
relatively free algebras of the Molien formula from Theorem 1.6.

{\bf Theorem 3.9.} [36, Theorem 7] {\sl Let $R$ be a PI-algebra,
let $G$ be a finite subgroup of $GL_m$ and let
$\xi_1(g),\ldots,\xi_m(g)$ be the eigenvalues of $g \in G$. Then
the Hilbert series of the algebra of invariants of $G$ acting on
the relatively free algebra $F_m(R)$ is}
$$
{\rm Hilb}(F_m(R)^G,t) = {1\over \vert G\vert}\sum_{g\in G}{\rm
Hilb}(F_m(R),\xi_1(g)t,\ldots,\xi_m(g)t).
$$

In the special case of commutative algebras when $F_m(R) =
F[x_1,\ldots,x_m]$ we obtain that
$$
\prod_{i=1}^m{1\over 1 - \xi_i(g)t} = {1\over {\rm det}(1 - gt)}
$$
which gives immediately the Molien formula from Theorem 1.6. In the
other extremal case $F_m(R) = K\langle X_m\rangle$, i.e. when
$T_m(R) = \{0\}$ we obtain that
$$
{\rm Hilb}(K\langle X_m\rangle,\xi_1(g)t,\ldots,\xi_m(g)t) =
{1\over 1 - (\xi_1(g) +\ldots +\xi_m(g))t} =
{1\over 1 - {\rm tr}(g)t}
$$
and this gives Theorem 2.5. Tracing the proof of Theorem 3.9 it is
easy to see that the result holds in the more general setup for the
algebra of invariants $(K\langle X_m\rangle/I)^G$ where the ideal
$I$ of $K\langle X_m\rangle$ is stable under the $GL_m$-action on
$K\langle X_m\rangle$.
\medskip

\centerline{\bf 4. Groups Acting on Free Lie Algebras}
\medskip

Recall that the (nonassociative) algebra $R$ is called a {\bf Lie
algebra} if it satisfies the {\bf Jacobi identity}
$$
(x_1x_2)x_3 + (x_2x_3)x_1 + (x_3x_1)x_2 = 0
$$
and the {\bf anticommutative identity}
$$
x_1^2 = 0\, (\hbox{\rm or equivalently}\, x_1x_2 = -x_2x_1).
$$
A typical example of a Lie algebra is the algebra $R^{(-)}$
associated with an associative algebra $R$ with new multiplication
$[r_1,r_2] = r_1r_2 - r_2r_1$, $r_1,r_2 \in R$. The famous
Poincar\'e-Birkhoff-Witt theorem states that any Lie algebra is
isomorphic to a subalgebra of $R^{(-)}$ for some associative
algebra $R$. In particular by the Witt theorem the {\bf free Lie
algebra} $L_m = L(X_m)$ of rank $m$ is isomorphic to the Lie
subalgebra generated by $x_1,\ldots,x_m$ in the free associative
algebra $K\langle X_m\rangle$. For our purposes it is convenient to
consider $L(X_m)$ as a Lie subalgebra of $K\langle X_m\rangle$ and
to use square brackets for the multiplication in the Lie algebras.
As a vector space $L(X_m)$ inherits the (multi)grading of $K\langle
X_m\rangle$ and its Hilbert series is given by the Witt formula
$$
{\rm Hilb}(L_m,t) = \sum_{n=1}^{\infty}\sum_{d\vert
n}\mu(d)m^{n/d}{t^n\over n},
$$
where $\mu(d)$ is the M\"obius function. By the theorem of Shirshov
any subalgebra of a free Lie algebra is also free.

Repeating verbatim the arguments for associative algebras one
introduces polynomial identities, relatively free algebras,
invariants etc. for Lie algebras. For a background on polynomial
identities for Lie algebras we refer to [5]. The main difference
between the associative and Lie PI-algebras is that the existence
of a polynomial identity is a very strong restriction on the
associative algebra and does not affect so much the structure
and combinatorial properties of the Lie algebra. From this point of
view the Lie PI-algebras are closer to the groups satisfying
identical relations. As a result methods from the well developed
combinatorial group theory are transferred first to Lie algebras
and then to associative algebras. In this section we summarize some
results on invariants of finite linear groups acting on free and
relatively free Lie algebras. Since the (relatively) free Lie
algebras are in an intermediate position between (relatively) free
associative algebras and (relatively) free groups one should expect
that the properties of the invariants in the Lie case are also
somewhere between the invariants in the associative case and the
subgroups of fixed points under the action of a finite group of
automorphisms of the (relatively) free group.

By the Engel theorem if $R$ is a finite dimensional Lie algebra
satisfying the Engel condition, then $R$ is nilpotent, i.e. there
exists a positive integer $c$ such that $[x_1,\ldots,x_c] = 0$ is a
polynomial identity for $R$. By the result of Kostrikin [57] the
finitely generated Lie algebras satisfying the Engel polynomial
identity $x_2({\rm ad} x_1)^n = 0$ for some $n$ are also nilpotent.
Recently Zel'manov [97] has established that the Engel identity
implies the nilpotency of the Lie algebra without the restriction
for finite generating. It is easy to see that if $f(x_1,\ldots,x_m)
= 0$ is a polynomial identity for a Lie algebra $R$ such that
$f(x_1,\ldots,x_m) \in L(X_m)$ does not belong to the second member
$L''_m = [[L_m,L_m],[L_m,L_m]]$ of the derivation series of $L_m$
then $f(x_1,\ldots,x_m) = 0$ implies the Engel identity. Therefore
for a T-ideal $T_m(R)$ of $L_m$ either $T_m(R) \subseteq L''_m$ or
the relatively free Lie algebra $F_m(R)$ is nilpotent (and hence
finite dimensional). If $T_m(R_1) \subseteq T_m(R_2)$ are two
T-ideals in $L_m$ there is a canonical homomorphism
$$
\pi: F_m(R_1) \cong L_m/T_m(R_1) \longrightarrow L_m/T_m(R_2) \cong
F_m(R_2).
$$
If $G$ is a finite subgroup of $GL_m$ then
$\pi(F_m(R_1)^G)\subseteq F_m(R_2)^G$. It follows easily from the
Maschke theorem that every invariant of $F_m(R_2)$ can be lifted to
an invariant of $F_m(R_1)$, i.e. $\pi(F_m(R_1)^G) = F_m(R_2)^G$. As
a consequence if the Lie algebra $R$ is not nilpotent and $G$ is a
finite linear group acting on $L_m$, then the properties of the
invariants of $G$ in $F_m(R)$ should be similar to those in the
free Lie algebra $L_m$ and in the free {\bf metabelian Lie algebra}
$L_m/L''_m$. The structure of $L_m/L''_m$ is well known. It has a
basis consisting of all commutators
$$
[x_{i_1},x_{i_2},\ldots,x_{i_n}],\, m \geq i_1 > i_2 \leq ... \leq
i_n \leq m,\, n = 1,2,\ldots\, .
$$
It is easy to calculate the Hilbert series of $L_m/L''_m$ (see e.g.
[24])
$$
{\rm Hilb}(L_m/L''_m,t_1,\ldots,t_m) = 1 + (t_1 + \ldots + t_m) +
(t_1 + \ldots + t_m - 1)\prod_{i=1}^m{1\over 1 - t_i}.
$$

Bryant [12] has proved that for any nontrivial finite subgroup $G$
of $GL_m$ the algebra of invariants $L_m^G$ is not finitely
generated. The partial case for $G$ cyclic has been also
established by Belov [6]. We should compare this result with the
theorem of Dyar and Scott [29] that for any finite group $G$ of
automorphisms of the free group $F_m$ the group of fixed points
$F_m^G = \{f \in F_m \mid g(f) = f,\, g \in G\}$ is finitely
generated. In the same paper [12] Bryant has also given an analogue
of the Witt formula for the subalgebras of $L_m$ generated by a
finite system of homogeneous elements. The approach of Bryant is
characteristic free and his results hold for free Lie algebras over
a principal ideal domain $K$. The author [24] has handled the other
extreme case $(L_m/L''_m)^G$ and as a consequence the case
$F_m(R)^G$ for $F_m(R)$ being non-nilpotent; see [24] for further
generalizations.

{\bf Theorem 4.1.} [12, 24] {\sl Let $G$ be a nontrivial finite
subgroup of $GL_m$, $m \geq 2$, and let $R$ be a Lie algebra. The
algebra of invariants $F_m(R)^G$ is finitely generated if and only
if $R$ is nilpotent.}

Another result on groups acting on free metabelian algebras is due
to Bryant [13] and may be useful in noncommutative invariant
theory. Let $G$ be a finite group. The $G$-module $W$ is {\bf free}
if it is a direct sum of $G$-submodules $W_i$, $i \in I$, and every
$W_i$ is isomorphic to the {\bf group algebra} $KG$ considered as a
left $G$-module. This means that for every $i \in I$ there exists
an element $w_i \in W_i$ such that as a vector space $W_i$ has a
basis $\{g(w_i) \mid g \in G\}$.

{\bf Theorem 4.2.} [13] {\sl Let $G$ be a finite subgroup of $GL_m$
and let $c$ be the number of elements of $G$ which act by scalar
multiplication on the free metabelian Lie algebra $M \cong
L_m/L''_m$. Let us denote by $W_n$, $n = 1,2,\ldots$, a maximal
free submodule of the $G$-module $M^n/M^{n+c}$ where $M^n$ is the
$n$th power of $M$, i.e. the vector space spanned on all
commutators of length $\geq n$. Then
$$
\lim_{n\to\infty}{{\rm dim}W_n\over{\rm dim}(M^n/M^{n+c})} = 1,
$$
i.e. for big $n$ the $G$-module $M^n/M^{n+c}$ is very close to a
free module.}
\medskip

\centerline{\bf 5. Generic Matrix Algebras}
\medskip

We start this section with a realization of the relatively free
algebra $F_m(M_k(K)) = K\langle X_m\rangle/T_m(M_k(K))$, where
$M_k(K)$ is the $k \times k$ matrix algebra with entries from the
base field $K$. The construction due to Procesi [2] has important
applications to structure theory of rings, invariant theory,
theory of division algebras, etc. (see e.g. [50] and [82]). Till the
end of the section we fix the positive integers $k$ for the size
of the matrices and $m$ for the number of generators. Let
$$
\Omega = K[x_{pq}^{(i)} \mid p,q = 1,\ldots,k,\, i = 1,\ldots,m]
$$
be the polynomial algebra in $k^2m$ commuting variables.

{\bf Definition 5.1.} {\sl The $k \times k$ matrices
$$
x_i = (x_{pq}^{(i)}) = \sum_{p=1}^k\sum_{q=1}^kx_{pq}^{(i)}e_{pq}
\in M_k(\Omega),\, i = 1,\ldots,m,
$$
where $e_{pq}$ are the matrix units, are called
{\bf generic matrices} and generate the {\bf generic matrix
algebra} $R_{km}$. We denote by $C_{km}$ the subalgebra of
$\Omega$ generated by all the traces of elements of $R_{km}$. The
{\bf generic trace algebra} $T_{km}$ is the subalgebra of
$M_k(\Omega)$ generated by $R_{km}$ and $C_{km}$ (where we
identify the traces with scalar matrices).}

The algebra $T_{km}$ is multigraded in a natural way counting the
appearings of $x_i$ in the monomials $x_{i_1} \ldots x_{i_n}$ and
their traces ${\rm tr}(x_{i_1} \ldots x_{i_n})$. For example
$$
{\rm tr}(x_2x_1x_3x_2){\rm tr}(x_1)x_2x_1x_2
$$
is homogeneous of multidegree $(3,4,1,0,\ldots,0)$.

{\bf Proposition 5.2.} {\sl The kernel of the canonical
epimorphism
$$
\pi: K\langle X_m\rangle\longrightarrow R_{km}
$$
defined by $\pi(x_i) = (x_{pq}^{(i)}) \in R_{km}$, $i =
1,\ldots,m$, is equal to $T_m(M_k(K))$, i.e. the generic $k \times
k$ matrix algebra $R_{km}$ is isomorphic to the relatively free
algebra $F_m(M_k(K))$.}

The {\bf adjoint action} of the group $GL_k$ on the generic
matrices $x_i \in R_{km}$ is defined as the action by conjugation
$$
g\ast x_i = gx_ig^{-1} = (f_{pq}^{(i)}),
$$
$g \in GL_k$, $f_{pq}^{(i)} \in \Omega$, $p,q = 1,\ldots,k$, $i =
1,\ldots,m$. Since the centre of $GL_k$ acts trivially on $x_i$ we
obtain an induced action of the projective special linear group
$PSL_k$ on
$$
{\rm span}\{x_{pq}^{(i)} \mid p,q = 1,\ldots,k,\, i = 1,\ldots,m\}
$$
defined by
$$
g:x_{pq}^{(i)}\longrightarrow f_{pq}^{(i)},\, g \in PSL_k, p,q =
1,\ldots,k,\, i = 1,\ldots,m,
$$
and therefore also an action of $PSL_k$ on the polynomial algebra
$\Omega$. The main purpose of this section is to survey some
properties of the algebra of invariants $\Omega^{PSL_k}$ which is
called the {\bf algebra of} ({\bf simultaneous}) {\bf invariants}
of $m$ generic $k \times k$ matrices. For more detailed exposition
of the results before 1986 we refer to the survey by Formanek [39]
(see also his book [40]). For other properties of the algebra
of the matrix invarinats interesting from the point of view of
commutative algebra, algebraic geometry and representation theory
of algebras see also the paper by Van den Bergh [93].

Since
$$
{\rm tr}(x_{i_1} \ldots x_{i_n}) = {\rm tr}(gx_{i_1}g^{-1}\ldots
gx_{i_n}g^{-1})
$$
for any monomial $x_{i_1} \ldots x_{i_n}\in R_{km}$ and any $g \in
GL_k$ we obtain that
$$
{\rm tr}(x_{i_1} \ldots x_{i_n}) \in \Omega^{PSL_k},
$$
i.e. $C_{km} \subseteq \Omega^{PSL_k}$. The first fundamental
theorem for the matrix invariants states that these two algebras
coincide.

{\bf Theorem 5.3.} (First Fundamental Theorem of Matrix
Invariants, [44, Theorem 16.2], [74, Theorem 1.3], [85, Theorem
1]) {\sl The algebra of simultaneous invariants $\Omega^{PSL_k}$
of $m$ generic $k \times k$ matrices is generated by all the
traces ${\rm tr}(x_{i_1}\ldots x_{i_n})$ of monomials from the
generic trace ring $T_{km}$, i.e. $C_{km} = \Omega^{PSL_k}$.}

Recall that if $G$ is a group and $V$ and $W$ are $G$-modules then
$G$ acts {\bf diagonally} on the tensor product $V \otimes W$ by
$$
g(v \otimes w) = g(v) \otimes g(w),\, g \in G, v \in V, w \in W.
$$
In particular if $G$ acts on an algebra $R$ and on the $k \times
k$ matrix algebra $M_k(K)$ we define the diagonal action on the
matrix algebra $M_k(R) \cong R \otimes_K M_k(K)$ with entries from
$R$ via
$$
g(\sum_{i=1}^k\sum_{j=1}^kr_{ij}e_{ij}) =
\sum_{i=1}^k\sum_{j=1}^kg(r_{ij})g(e_{ij}),
\, g \in G, r_{ij} \in R.
$$
The next result describes the generic trace ring as an algebra of
matrix concominants.

{\bf Theorem 5.4.} ([74, Section 2], [35, Theorem 10]) {\sl Let
$GL_k$ act on $\Omega$ as described above and on $M_k(K)$ by $g(a)
= g^{-1}ag$, where $a \in M_k(K)$ and $g \in GL_k$. (Note that
$g^{-1}$ and $g$ are reversed in comparison to the action of
$GL_k$ on $\Omega$ induced by the action on the generic matrices
$x_i = (x_{pq}^{(i)})$.) This induces a diagonal action of $GL_k$
on $M_k(\Omega) \cong \Omega \otimes_K M_k(K)$ and the generic
trace algebra is the fixed algebra of this action.}

Recently Seelinger [83] has generalized Theorem 5.4 and has found
the invariants of the action of $GL_k$ on the tensor powers
$M_k(\Omega)^{\otimes n}$ of $M_k(\Omega)$. We refer to his paper
for additional information.

In virtue of Theorem 1.4 the algebra of invariants
$\Omega^{PSL_k}$ is finitely generated and the problem which
arises immediately is to find a bound for the degree of the
generators. For this purpose we recall the classical
Dubnov-\-Ivanov-\-Nagata-\-Higman Theorem [28, 69, 46] on the
class of nilpotency of nil algebras of bounded index.

{\bf Theorem 5.5.} {\sl Let $R$ be a non-unitary associative
algebra over the field $K$ of characteristic $0$ and let $R$
satisfy the polynomial identity $x^k = 0$. There exists an integer
$d = d(k)$ depending on $k$ only such that $R$ is nilpotent of
class $d$, i.e. $R$ satisfies the polynomial identity $x_1 \ldots
x_d = 0$.}

It turns out that there is a nice relation between the degree of
the generators of $C_{km}$ and the class of nilpotency $d(k)$ from
Theorem 5.5.

{\bf Theorem 5.6.} (Formanek [38, Theorem 6], Procesi [74,
Theorem 3.3], [75, Theorem 4.3], Razmyslov [80, Final remark])
{\sl For a fixed positive integer $k$ the algebra $C_{km}$ is
generated by the set of traces ${\rm tr}(x_{i_1} \ldots x_{i_n})$
of degree $n$ bounded by the class of nilpotency $d(k)$ of the nil
algebras of index $k$. For $m$ large enough the bound $n \leq d(k)$ is
exact.}

The best known bounds for $d(k)$
$$
{k(k+1)\over 2}\leq d(k) \leq k^2
$$
are due respectively to Kuz'min [59] and Razmyslov [80]. The only
exact values of $d(k)$ are known for $k \leq 3$ [27]:
$$
d(1) = 1, d(2) = 3, d(3) = 6.
$$
A minimal set of generators for the invariants of $2 \times 2$
matrices has been found by Siberskii [85] (see Theorem 5.10
below). Silvana Abeasis and Marilena Pittaluga [1] have
suggested an algorithm for finding a minimal set of generators for
$C_{km}$. Combining computers with calculations by hand they have
successfully aplied this algorithm to the case of $3 \times 3$
matrices.

The Second Fundamental Theorem or the problem for the description
of the defining relations of the algebra of $k \times k$ matrix
invariants (see Theorem 1.3) has two aspects. One of them is to
handle the case of the algebra $C_{km}$ for fixed $m$. The other
is to consider the multilinear invariants because every T-ideal
$T(R)$ of the free algebra of countable rank $K\langle
x_1,x_2,\ldots\rangle$ is generated as a T-ideal by its
multilinear elements and a similar assertion holds for the ideal
of relations of the invariants of $K[x_{pq}^{(i)} \mid i =
1,2,\ldots]$. The second aspect is completely solved by the
theorem of Helling-\-Razmyslov-\-Procesi [45, 80, 74]. For details
we refer to [39] and state a version of the result only. We need
some preliminary discussions.
For a fixed $n$ we write the permutations $\sigma \in S_n$ as
products of disjoint cycles
$$
\sigma = (i_1 \ldots i_p) \ldots (j_1 \ldots j_q),
$$
including also the 1-cycles, so that each integer $1,\ldots,n$
occures exactly once. We define the {\bf associated trace
function}
$$
{\rm tr}_{\sigma}(x_1 \ldots x_n) = {\rm tr}(x_{i_1} \ldots
x_{i_p}) \ldots {\rm tr}(x_{j_1} \ldots x_{j_q}),
$$
where $x_1,\ldots,x_n$ are generic $k \times k$ matrices. We also
assume that for $m \leq n$ the symmetric group $S_m$ acts on
$1,\ldots,m$ and leaves invariant $m + 1,\ldots,n$, i.e. $S_m$ is
canonically embedded into $S_n$.

{\bf Theorem 5.7.} (The Second Fundamental Theorem of Matrix
Invariants, [45, 80, 74]) {\sl Let
$$
f(x_1,\ldots,x_n) = \sum_{\sigma\in S_n}\alpha_{\sigma}{\rm
tr}_{\sigma}(x_1 \ldots x_n),\, \alpha_{\sigma} \in K,
$$
be a multilinear trace polynomial of degree $n$. Then $f = 0$ is
trace identity for the $k \times k$ matrix algebra, i.e.
$f(a_1,\ldots,a_n) = 0$ for all $a_1,\ldots,a_n \in M_k(K)$ if and
only if $\sum_{\sigma\in S_n}\alpha_{\sigma}\sigma$ belongs to the
two-sided ideal of the group algebra $KS_n$ generated by the
element $\sum_{\sigma\in S_{k+1}}({\rm sign} \sigma)\sigma$.}

The {\bf fundamental trace identity}
$$
\sum_{\sigma\in S_{k+1}}({\rm sign} \sigma){\rm
tr}_{\sigma}(x_1 \ldots x_n)=0
$$
is actually the linearization of the Cayley-Hamilton Theorem. We
illustrate this for $2 \times 2$ matrices. For any matrix $x \in
M_2(K)$ the Cayley-\-Hamilton Theorem states that
$$
x^2 - {\rm tr}(x)x + {\rm det}(x) = 0.
$$
Since the base field is algebraically closed the matrix $x$ is
conjugated with an upper triangular matrix
$$
\pmatrix{\xi_1&\eta\cr
0&\xi_2\cr}
$$
and ${\rm tr}(x) = \xi_1 + \xi_2$, ${\rm det}(x) = \xi_1\xi_2 =
{1\over 2} ({\rm tr}^2(x) - {\rm tr}(x^2))$ and we rewrite the
Caley-\-Hamilton Theorem as
$$
c(x) = x^2 - {\rm tr}(x)x + {1\over 2}({\rm tr}^2(x) - {\rm
tr}(x^2)) = 0.
$$
Now we linearize the identity $c(x) = 0$, i.e. consider the trace
identity $c(x_1 + x_2) - c(x_1) - c(x_2) = 0$ and obtain
$$
h(x_1,x_2) = x_1x_2 + x_2x_1 - {\rm tr}(x_1)x_2 - {\rm tr}(x_2)x_1
+ {\rm tr}(x_1){\rm tr}(x_2) - {\rm tr}(x_1x_2) = 0.
$$
Since the trace is a non-degenerated bilinear form on $M_2(K)$, the
vanishing of the polynomial $h(x_1,x_2)$ on $M_2(K)$ is equivalent
to the vanishing of ${\rm tr}(h(x_1,x_2)x_3)$ on all $2 \times 2$
matrices, i.e. ${\rm tr}(h(x_1,x_2)x_3)$ is equal to 0 in the
trace ring $T_{2m}$, $m \geq 3$. Direct calculations show that
$$
0 = {\rm tr}(h(x_1,x_2)x_3) = {\rm tr}(x_1x_2x_3) + {\rm
tr}(x_2x_1x_3) - {\rm tr}(x_1){\rm tr}(x_2x_3) -
$$
$$
- {\rm tr}(x_2){\rm tr}(x_1x_3) + {\rm tr}(x_1){\rm tr}(x_2){\rm
tr}(x_3) - {\rm tr}(x_1x_2){\rm tr}(x_3) =
$$
$$
  = {\rm tr}_{(123)}(x_1,x_2,x_3) + {\rm tr}_{(213)}(x_1,x_2,x_3)
- {\rm tr}_{(23)(1)}(x_1,x_2,x_3) -
$$
$$
-{\rm tr}_{(13)(2)}(x_1,x_2,x_3) + {\rm
tr}_{(1)(2)(3)}(x_1,x_2,x_3) - {\rm tr}_{(12)(3)}(x_1,x_2,x_3)=
$$
$$
= \sum_{\sigma\in S_3}({\rm sign} \sigma){\rm tr} (x_1,x_2,x_3).
$$

Considering the defining relations of the algebra $C_{km}$ for
fixed $m$ it is known that the transcendence degree of $C_{km}$ is
$$
{\rm trans.deg.}(C_{km}) = (m - 1)k^2 + 1,
$$
i.e. the algebra $C_{km}$ contains $(m - 1)k^2 + 1$ algebraically
independent elements and every $(m - 1)k^2 + 2$ elements satisfy a
non-trivial algebraic equation. For $2 \times 2$ matrices this
result sounds as follows.

{\bf Theorem 5.8.} [4, Theorem 6] {\sl The polynomials
$$
{\rm tr}(x_p), {\rm tr}(x_p^2), {\rm tr}(x_1x_p), {\rm
tr}(x_2x_q),\, p = 2,\ldots,m, q = 3,\ldots,m,
$$
form a maximal set of algebraically independent elements of
$F_{2m}$.}

Let $Q(C_{km})$ be the field of quotients of $C_{km}$. A well known
and important problem is  {\sl  whether  $Q(C_{km})$  is  a  purely
transcendental extension of the base field $K$, i.e. if there exist
$n  =  {\rm  trans.deg.}(C_{km})$  elements  $u_1,\ldots,u_n$  from
$Q(C_{km})$ such that $Q(C_{km})$ is isomorphic  to  the  field  of
rational functions $K(u_1,\ldots,u_n)$.} Procesi  ([72],  see  also
his  book  [73])  has  shown   that   $Q(C_{km})$   is   a   purely
transcendental extension of $Q(C_{k2})$, i.e. the problem should be
solved for $m = 2$ only. He has also handled  the  case  $k  =  2$.
Formanek [32, 33] has establshed the result for $k = 3, 4$. Another
approach to the rationality of $Q(C_{km})$  is  given  by  Van  den
Bergh [94] which gives an alternative proof for $k = 3$.  A  survey
on this topic can be found in [64].

{\bf Theorem 5.9.} [72, 32, 33] {\sl The field of quotients
$Q(C_{km})$ of the algebra of $k \times k$ matrix invariants is a
purely transcendental extension of the base field $K$ for $k =
2,3,4$.}

Explicit results for the generators and defining relations of
$C_{km}$ are known for $2 \times 2$ matrices only.

{\bf Theorem 5.10.} (i) (First Fundamental Theorem of $2 \times 2$
Matrix Invariants, Siberskii [85]) {\sl The polynomials
$$
{\rm tr}(x_p),\, p = 1,\ldots,m,\, {\rm tr}(x_px_q),\, 1 \leq p
\leq q \leq m,\, {\rm tr}(x_px_qx_r),\, 1 \leq p < q < r \leq m,
$$
form a minimal system of generators of $C_{2m}$.}

(ii) (Second Fundamental Theorem for the Algebra of
Invariants of Two and Three Generic $2 \times 2$ Matrices,
Siberskii [85], Formanek [35]) {\sl Let
$$
a_p = {\rm tr}(x_p), b_p = {\rm det}(x_p), c_p = {\rm tr}(x_qx_r),
d = {\rm tr}(x_1){\rm tr}(x_2){\rm tr}(x_3),
$$
$p,q,r = 1,2,3$, $p \not= q,r$, $q < r$.
The algebra of invariants $C_{22}$
of two generic $2 \times 2$ matrices is isomorphic to the
polynomial algebra $K[a_1,a_2,b_1,b_2,c_3]$. The algebra of
invariants of three generic matrices $C_{23}$ is generated by the
commuting elements $a_p,b_p,c_p,d$, $p = 1,2,3$, modulo the
principal ideal generated by the relation
$$
d^2 - (a_1c_1 + a_2c_2 + a_3c_3 - a_1a_2a_3)d +
$$
$$
+ b_1c_1^2 + b_2c_2^2 + b_3c_3^2 - a_1a_2b_3c_3 - a_3a_1b_2c_2 -
a_2a_3b_1c_1 +
$$
$$
+ a_1^2b_2b_3 + a_2^2b_3b_1 + a_3^2b_1b_2 - 4b_1b_2b_3 + c_1c_2c_3
= 0.
$$
Alternatively $C_{23}$ is a free module of rank $2$ (with basis
$1,d$) over the polynomial algebra $K[a_p,b_p,c_p \mid p =
1,2,3]$.}

The last statement of Theorem 5.10 (ii) has been generalized by
Teranishi for $k = 2$ any $m$ [91]. Van den Bergh [93] has
proved a stronger result that for any $k$ and $m$ the algebra of
invariants $C_{km}$ is Cohen-Macaulay, i.e. a free module over a
graded polynomial subalgebra of $C_{km}$.

Another approach to the defining relations of $T_{km}$ and
$C_{km}$ has been suggested by the author and Koshlukov [26] and
[25]. The idea has its sources from the theory of the identities
for representations of groups and algebras and works successfully
for $2 \times 2$ matrices.

{\bf Definition 5.11.} {\sl Let $N$ be a $GL_m$-module with basis
as a vector space $y_1,\ldots,y_p$. The ideal $I$ of the
polynomial algebra $K[y_1,\ldots,y_p]$ or of the free associative
algebra $K\langle y_1,\ldots,y_p\rangle$ is called a $GL$-{\bf
ideal} if it is invariant under the action of the general linear
group $GL_m$ induced on $K[y_1,\ldots,y_p]$ and $K\langle
y_1,\ldots,y_p\rangle$ by the action of $GL_m$ on $N$.}

The trace algebra $T_{km}$ and the algebra of matrix invariants
are $GL_m$-modules with the action of $GL_m$ induced by the action
of $GL_m$ on ${\rm span}\{x_1,\ldots,x_m\}$:
$$
g({\rm tr}(x_{h_1} \ldots x_{h_p})\ldots{\rm tr}(x_{i_1} \ldots
x_{i_q})x_{j_1} \ldots x_{j_r}) =
$$
$$
={\rm tr}(g(x_{h_1}) \ldots g(x_{h_p}))\ldots {\rm tr}(g(x_{i_1})
\ldots g(x_{i_q}))g(x_{j_1}) \ldots g(x_{j_r}),
$$
$g \in GL_m$, ${\rm tr}(x_{h_1} \ldots x_{h_p})\ldots {\rm
tr}(x_{i_1} \ldots x_{i_q})x_{j_1} \ldots x_{j_r} \in T_{km}$. If
$d = d(k)$ is the integer from Theorem 5.5 and
$$
N = {\rm span}\{{\rm tr}(x_{i_1} \ldots x_{i_n}) \mid i_p =
1,\ldots,m, p = 1,\ldots,n \leq d\}
$$
then Theorem 5.6 gives that $C_{km}$ and $T_{km}$ are generated by
$N$ and $N \oplus {\rm span}\{x_1,\ldots,x_m\}$, respectively,
both vector spaces are $GL_m$-modules. Hence it is natural
to search for a minimal system of generators and defining relations
of $T_{km}$ and $C_{km}$ in the language of $GL$-ideals. The
algorithm of Silvana Abeasis and Marilena Pittaluga [1] can be used
to find the structure of the minimal $GL_m$-submodule of $N$
generating $C_{km}$. The description in [1] is in the language of
representation theory of $S_n$ but it is easy to translate the
obtained results in the language of $GL_m$ as in [19, 7].

{\bf Definition 5.12.} (Razmyslov [77, 78]) {\sl Let $R$ be an
associative algebra and let $H$ be a Lie subalgebra of $R$ which
generates $R$ as an associative algebra. The element
$f(x_1,\ldots,x_m)$ from the free associative algebra $K\langle
x_1,\ldots,x_m\rangle$ is called a {\bf weak polynomial identity
for the pair} $(R,H)$ if $f(h_1,\ldots,h_m) = 0$ in $R$ for every
$h_1,\ldots,h_m \in H$.}

The following theorems give the complete description of the
defining relations of $T_{2m}$ and $C_{2m}$.

{\bf Theorem 5.13.} {\sl Let $m \geq 2$ and let
$$
y_p = x_p - {1\over 2}{\rm tr}(x_p) =
\pmatrix{y_{11}^{(p)}&y_{12}^{(p)}\cr
y_{21}^{(p)}&-y_{11}^{(p)}\cr},
$$
where $y_{11}^{(p)} = {1\over 2}(x_{11}^{(p)} - x_{22}^{(p)})$,
$y_{12}^{(p)} = x_{12}^{(p)}$, $y_{21}^{(p)} = x_{21}^{(p)}$, $p =
1,\ldots,m$, be the {\bf generic traceless $2 \times 2$ matrices}.}

(i) (Procesi [76]) {\sl The generic trace algebra $T_{2m}$ is
isomorphic to the tensor product of $K$-algebras
$$
K[{\rm tr}(x_1),\ldots,{\rm tr}(x_m)] \otimes_K R'_{2m},
$$
where $R'_{2m}$ is the associative algebra generated by the
generic traceless matrices $y_1,\ldots,y_m$.}

(ii) (Razmyslov [77]) {\sl The algebra of the generic traceless
matrices $R'_{2m}$ is isomorphic to the factor-algebra $K\langle
y_1,\ldots,y_m\rangle/T_m(M_2(K),sl_2(K))$ of the free associative
algebra $K\langle y_1,\ldots,y_m\rangle$, where $sl_2(K)$ is the
Lie algebra of all traceless $2 \times 2$ matrices and
$T_m(M_2(K),sl_2(K))$ is the ideal of all weak polynomial
identities in $K\langle y_1,\ldots,y_m\rangle$ for the pair
$(M_2(K),sl_2(K))$. The ideal $T_m(M_2(K),sl_2(K))$ is generated
as a weak T-ideal by the weak polynomial identity $[x_1^2,x_2] =
0$, i.e. $T_m(M_2(K),sl_2(K))$ is the minimal ideal of weak
polynomial identities containing the element $[x_1^2,x_2]$.}

(iii) (Drensky and Koshlukov [26, Corollary 2]) {\sl Let $GL_m$
act canonically on the free associative algebra $K\langle
y_1,\ldots,y_m\rangle$. As a $GL$-ideal $T_m(M_2(K),sl_2(K))$ is
generated by the weak polynomial identities $[x_1^2,x_2]$ and
$s_4(x_1,x_2,x_3,x_4)$, the second relation appears for $m \geq 4$
only. Hence the algebra of $2 \times 2$ generic traceless matrices
has a uniform set of defining relations for any $m \geq 2$.}

{\bf Theorem 5.14.} (Second Fundamental Theorem of $2 \times 2$
Matrix Invariants, Drensky [25]) {\sl Let $m \geq 2$, let
$y_1,\ldots,y_m$ be the $2 \times 2$ generic traceless matrices
and let $GL_m$ act on the vector space spanned on $a_p = {\rm
tr}(x_p)$, $b_{pq} = {\rm tr}(y_py_q)$, $c_{pqr} = {\rm
tr}(s_3(y_p,y_q,y_r))$, $p,q,r = 1,\ldots,m$, as indicated in the
text below Definition 5.11. The algebra of invariants $C_{2m}$
is isomorphic to the polynomial algebra
$$
K[a_d, b_{ij}, c_{pqr} \mid d,i,j,p,q,r = 1,\ldots,m, i \leq j, p
< q < r]
$$
modulo the $GL$-ideal generated by the relations (some of them may
miss for small $m$)
$$
u_1 = 2c_{123}^2 - 3\sum_{\pi\in S_3}({\rm sign}
\pi)b_{1\pi(1)}b_{2\pi(2)}b_{3\pi(3)},\, m \geq 3,
$$
$$
u_2 = -b_{11}c_{234} + b_{12}c_{134} - b_{13}c_{124} +
b_{14}c_{123},\, m \geq 4,
$$
$$
u_3 = \sum({\rm sign}
\pi)c_{\pi(1)\pi(2)\pi(3)}c_{1\pi(4)\pi(5)},\, m \geq 5,
$$
where in the third relation the summation runs over all
permutations $\pi \in S_5$ such that $\pi(1) < \pi(2) < \pi(3)$,
$1 < \pi(4) < \pi(5)$, i.e. $C_{2m}$ is generated by ${\rm
tr}(x_p)$, ${\rm tr}(y_py_q)$ and ${\rm tr}(s_3(y_p,y_q,y_r))$
with $GL$-defining relations
$$
u_1 = 2{\rm tr}^2(s_3(y_1,y_2,y_3)) -
3\sum_{\pi\in S_3}({\rm sign} \pi){\rm tr}(y_1,y_{\pi(1)}){\rm
tr}(y_2y_{\pi(2)}){\rm tr}(y_3y_{\pi(3)}) = 0,
$$
$$u_2 = \sum_{p=1}^4(-1)^p{\rm tr}(s_3(y_1,\ldots,\hat
y_p,\ldots,y_4){\rm tr}(y_1y_p) = 0,
$$
( $\hat{}$ means that the corresponding $y_p$ does not appear),}
$$
u_3 = \sum_{\pi\in S_5}({\rm sign} \pi){\rm
tr}(s_3(y_{\pi(1)},y_{\pi(2)},y_{\pi(3)})){\rm
tr}(s_3(y_1,y_{\pi(4)},y_{\pi(5)}))=0.
$$

Now we survey some results devoted to the Hilbert series of
$T_{km}$ and $C_{km}$. In virtue of Theorems 1.4 and 5.6 $C_{km}$
is a commutative algebra generated by a finite system of
multihomogeneous elements. A general result of commutative algebra
gives that the Hilbert series of $C_{km}$ is a rational function
of $t_1,\ldots,t_m$. The algebra $T_{km}$ is a finitely generated
$C_{km}$-module and by the same result its Hilbert series is also
rational. The Weyl analogue for compact groups of the Molien
formula 1.6 [96] gives formulae for the Hilbert series of $C_{km}$
and $T_{km}$ as (multiple) integrals. Recall that the complex
unitary group $U(k,{\bf C})$ consists of all invertible $k \times
k$ matrices $g$ with complex entries such that $g^{-1} = ^t\bar
g$, where $^t\bar g$ is the transposed matrix of the complex
conjugated of $g$.

{\bf Theorem 5.15.} (see e.g. [35, p. 204]) {\sl Let $U(k,{\bf
C})$ be the complex unitary group with normalized Haar measure
$\mu$. Then
$$
{\rm Hilb}(C_{km},t_1,\ldots,t_m) = \int_{U(k,{\bf C})}
\left(\prod_{i=1}^m{1\over{\rm det}(1 - t\phi^g)}\right)d\mu(g),
$$
$$
{\rm Hilb}(T_{km},t_1,\ldots,t_m) = \int_{U(k,{\bf C})}
{\rm tr}(\phi^g)
\left(\prod_{i=1}^m{1\over{\rm det}(1 - t\phi^g)}\right)d\mu(g),
$$
where for each $g \in U(k,{\bf C})$, $\phi^g \in U(k^2,{\bf C})$ is
a matrix giving the action of $g$ by conjugation on $M_k({\bf
C})$.}

Unfortunately the evaluation of the multiple integrals from
Theorem 5.15 is difficult. Another description of the Hilbert
series of $C_{km}$ and $T_{km}$ in the language of representation
theory of $GL_m$ is given by Formanek [35]; see also [34] for
references how to translate the results of Razmyslov-\-Procesi
[80, 74] in order to obtain another expression of the Hilbert
series of $C_{km}$ and $R_{km}$. In particular Formanek has
succeeded to evaluate explicitly $C_{2m}$ and $T_{2m}$ (see also
[76] and the book of Le Bruyn [63] for other approaches to
$C_{2m}$ and $T_{2m}$). It is known that the $GL_m$-module
structure of $C_{2m}$ and $T_{2m}$ is completely determined by
their Hilbert series for $m = 4$ only and we give the result in
this form.

{\bf Theorem 5.16.} ([35, Theorem 12], see also [76] and [63])
$$
{\rm Hilb}(C_{24},t_1,t_2,t_3,t_4) =
{1 + \sum t_pt_qt_r - t_1t_2t_3t_4\sum t_p - (t_1t_2t_3t_4)^2
\over\prod(1 - t_p)(1 - t_p^2)\prod(1 - t_pt_q)},
$$
$$
{\rm Hilb}(T_{24},t_1,t_2,t_3,t_4) =
{1 - t_1t_2t_3t_4
\over\prod(1 - t_p)^2\prod(1 - t_pt_q)},
$$
{\sl where the sums and products are over $1 \leq p \leq 4$, $1
\leq p < q \leq 4$, $1 \leq p < q < r \leq 4$, respectively.}

Teranishi [89, 90] has applied the Cauchy integral formula to Theorem
5.15 and has obtained the explicit form of the Hilbert series of
$C_{k2}$ and $T_{k2}$ for $k = 3$ and $4$.

Since the generic matrix algebra $R_{km}$ is a subalgebra of the
generic trace algebra $T_{km}$ one may uses the information for
$T_{km}$ and $C_{km}$ in order to describe $R_{km}$ and its centre
$R_{km} \cap C_{km}$. We are not going to discuss this topic in
detail here. We shall only mention the famous and important
problem of Kaplansky from 1956 [51] {\sl whether the centre of
$R_{km}$ is non-trivial for any $k \geq 3$ and $m \geq 2$.} It has
been solved positively by Formanek [39] and Razmyslov [78]
and has forced a revision of the structure theory of PI-algebras
(see [50, 82]). We refer to [34, 40, 22] for a survey on
the results on central polynomials. Concerning the Hilbert series
of $R_{km}$ and $R_{km} \cap C_{km}$, Formanek [35] has proved
that the coefficients of the Hilbert series of the graded vector
spaces $T_{km}/R_{km}$ and $C_{km}/(R_{km} \cap C_{km})$ are
``small'', i.e. the dimensions of the homogeneous components of
degree $n$ of $T_{km}$ and $R_{km}$ (and $C_{km}$ and $R_{km} \cap
C_{km}$) are asymptotically equal for $n\longrightarrow\infty$ and
for $k$ and $m$ fixed.

Since $C_{km}$ and $R_{km}$ have additional nice algebraic
properties their Hilbert series should be nicely looking. Van den
Bergh [95] has used some ideas of Stanley [88] to reduce the
determination of the Hilbert series of $C_{km}$ and $R_{km}$ to a
problem about flows in a certain graph and has obtained an
important consequence for the denominators of the rational
functions giving the explicit form of the series.

{\bf Theorem 5.17.} [95] {\sl The Hilbert series of $C_{km}$ and
$T_{km}$ can be expressed as rational functions whose denominators
are products of terms $(1 - u)$ where $u$ is a monomial in
$t_1,\ldots,t_m$ of degree $\leq k$.}

Another confirmation for the nice properties of $C_{km}$ and
$T_{km}$ is that their Hilbert series satisfy a functional
equation.

{\bf Theorem 5.18.} (Formanek [37], Teranishi [89, 91], see also Le
Bruyn [62] for the case $k = 2$) {\sl Let $H(t_1,\ldots,t_m)$ be
the Hilbert series of either $C_{km}$ or $T_{km}$, where $m \geq
2$ for $k\geq 3$ and $m > 2$ for $k = 2$. Then $H(t_1,\ldots,t_m)$
satisfies the functional equation
$$
H(t_1^{-1},\ldots,t_m^{-1}) = (-1)^d(t_1 \ldots
t_m)^{k^2}H(t_1,\ldots,t_m),
$$
where $d$ is the Krull dimension of $C_{km}$ (or respectively of
$T_{km}$).} (It is known that $d = {\rm trans.deg.}(C_{km}) = (m -
1)k^2 + 1$.)

Obviously every subgroup of $GL_m$ also acts on the vector space
${\rm span}\{x_{pq}^{(i)} \mid p,q = 1,\ldots,k, i = 1,\ldots,m\}$
and it is worthy to consider the algebra of invariants $\Omega^G$
of some other classical groups $G$ on $\Omega =  K[x_{pq}^{(i)}
\mid p,q = 1,\ldots,k, i = 1,\ldots,m]$. Such investigations are
done for example for $Sp(k)$, $SO(k)$, $O(k)$ and involve generic
matrix algebras with involution. We refer e.g. to the papers of
Procesi [74], Berele [8, 9, 10], H. Aslaksen, Eng-Chye Tan, Chen-Bo
Zhu (see [4] and the references there) for details.

\medskip
\centerline{\bf References}
\par
\item{1.}
S. Abeasis, M. Pittaluga,
On a minimal set of generators for the invariants of $3 \times 3$
matrices,
Commun. Algebra 17 (1989), 487 -- 499.
\par
\item{2.}
S.A. Amitsur, C. Procesi,
Jacobson rings and Hilbert algebras,
Ann. Mat. Pura Appl. 71 (1966), 67 -- 72.
\par
\item{3.}
A.Z. Anan'in,
Locally finitely approximable and locally representable varieties of algebras (Russian),
Algebra i Logika 16 (1977), 3 -- 23.
Translation: Algebra and Logic 16 (1977), 1 -- 16.
\par
\item{4.}
H. Aslaksen, Eng-Chye Tan, Chen-Bo Zhu,
Generators and relations of invariants of $2 \times 2$ matrices,
Comm. Algebra 22 (1994), 1821 -- 1832.
\par
\item{5.}
Yu.A. Bakhturin,
Identical Relations in Lie Algebras (Russian),
Nauka, Moscow, 1985.
Translation: VNU Science Press, Utrecht, 1987.
\par
\item{6.}
A.I. Belov,
Periodic automorphisms of free Lie algebras and their fixed points,
preprint.
\par
\item{7.}
A. Berele,
Homogeneous polynomial identities,
Israel J. Math. 42 (1982), 258 -- 272.
\par
\item{8.}
A. Berele,
Matrices with involution and nvariant theory,
J. Algebra 135 (1990), 139 -- 164.
\par
\item{9.}
A. Berele,
Flows on graphs with two-headed and two-tailed edges,
preprint.
\par
\item{10.}
A. Berele,
Trace identities and ${\bf Z}/2{\bf Z}$-graded invariants,
preprint.
\par
\item{11.}
A. Braun,
The nilpotency of the radical in a finitely generated P.I. ring,
J. Algebra 89 (1984), 375 -- 396.
\par
\item{12.}
R.M. Bryant,
On the fixed points of a finite group acting on a free Lie algebra,
J. London Math. Soc. (2) 43 (1991), 215 -- 224.
\par
\item{13.}
R.M. Bryant,
Symmetric powers of representations of finite groups,
J. Algebra 154 (1993), 416 -- 436.
\par
\item{14.}
L. Carini,
The Poincar\'e series related to the Grassmann algebra,
Lin. Multilin. Algebra 27 (1990), 199 -- 205.
\par
\item{15.}
C. Chevalley,
Invariants of finite groups generated by reflections,
Amer. J. Math. 67 (1955), 778 -- 782.
\par
\item{16.}
W. Dicks, E. Formanek,
Poicar\'e series and a problem of S. Montgomery,
Lin. Multilin. Algebra 12 (1982), 21 -- 30.
\par
\item{17.}
J.A. Dieudonn\'e, J.B. Carrell,
Invariant Theory, Old and New,
Academic Press, New York - London, 1971.
\par
\item{18.} M. Domokos,
Relatively free invariant algebras of finite reflection groups,
Trans. Am. Math. Soc. 348 (1996), No. 6, 2217 -- 2234.
\par
\item{19.}
V. Drensky,
Representations of the symmetric group and varieties of linear algebras (Russian),
Mat.Sb. 115 (1981), 98 -- 115.
Translation: Math. USSR Sb. 43 (1981), 85 -- 101.
\par
\item{20.}
V. Drensky,
Codimensions of T-ideals and Hilbert series of relatively free
algebras,
J. Algebra 91 (1984), 1 -- 17.
\par
\item{21.}
V. Drensky,
Polynomial identities of finite dimensional algebras,
Serdica 12 (1986), 209 -- 216.
\par
\item{22.}
V. Drensky,
Methods of commutative algebra in non-commutative ring theory,
Atti Accad. Peloritana dei Pericolanti, Classe I di Sci. Fiz., Mat.
e Nat. 70 (1992), Parte Prima, 239 -- 260.
\par
\item{23.}
V. Drensky,
Finite generation of invariants of finite linear groups on
relatively free algebras,
Lin. and Multilin. Algebra 35 (1993), 1 -- 10.
\par
\item{24.}
V. Drensky,
Fixed algebras of residually nilpotent Lie algebras,
Proc. Amer. Math. Soc. 120 (1994), 1021 -- 1028.
\par
\item{25.}
V. Drensky,
Defining relations for the algebra of invariants of $2 \times 2$ matrices,
Algebras and Representation Theory 6 (2003), No. 2, 193 -- 214.
\par
\item{26.}
V. Drensky, P. Koshlukov,
Weak polynomial identities for a vector space with a symmetric
bilinear form,
Math. and Education in Math.
Proc. of 16-th Spring Conf. of the Union of Bulg. Mathematicians,
Sofia, Publ. House of the Bulg. Acad. of Sci., 1987, 213 -- 219.
\par
\item{27.}
J. Dubnov,
Sur une g\'en\'eralisation de l'\'equation de Hamilton - Cayley et
sur les invariants simultan\'es de plusieurs affineurs,
Proc. Seminar on Vector and Tensor Analysis, Mechanics Research
Inst., Moscow State Univ. 2/3 (1935), 351 -- 367
(see also Zbl. f\"ur Math. 12 (1935), p. 176).
\par
\item{28.}
J. Dubnov, V. Ivanov,
Sur l'abaissement du degr\'e des polyn\^omes en affineurs,
C.R. (Doklady) Acad. Sci. USSR 41 (1943), 96 -- 98
(see also MR 6 (1945), p.113, Zbl. f\"ur Math. 60 (1957), p. 33).
\par
\item{29.}
J.L. Dyer, G.P. Scott,
Periodic automorphisms of free groups,
Commun. Algebra 3 (1975), 195 -- 201.
\par
\item{30.}
J.W. Fisher, S. Montgomery,
Invariants of finite cyclic groups
acting on generic matrices,
J. Algebra 99 (1986), 430 -- 437.
\par
\item{31.}
E. Formanek,
Central polynomials for matrix rings,
J. Algebra 23 (1972), 129 -- 132.
\par
\item{32.}
E. Formanek,
The center of the ring of $3 \times 3$ generic matrices,
Lin. Multilin. Algebra 7 (1979), 203 -- 212.
\par
\item{33.}
E. Formanek,
The center of the ring of $4 \times 4$ generic matrices,
J. Algebra 62 (1980), 304 -- 319.
\par
\item{34.}
E. Formanek,
The polynomial identities of matrices,
Contemp.Math. 13 (1982), 41 -- 79.
\par
\item{35.}
E. Formanek,
Invariants and the ring of generic matrices,
J. Algebra 89 (1984), 178 -- 223.
\par
\item{36.}
E. Formanek,
Noncommutative invariant theory,
Contemp. Math. 43 (1985), 87--119.
\par
\item{37.}
E. Formanek,
Functional equations for character series associated with
$n \times n$ matrices,
Trans. Amer. Math. Soc. 294 (1986), 647 -- 663.
\par
\item{38.}
E. Formanek,
Generating the ring of matrix invariants,
Lect. Notes in Math. 1195, Springer Verlag, Berlin - Heidelberg -
New York, 1986, 73 -- 82.
\par
\item{39.}
E. Formanek,
The invariants of $n \times n$ matrices,
Lect. Notes in Math. 1278, Springer Verlag, Berlin - Heidelberg -
New York - London - Paris - Tokyo, 1987, 18 -- 43.
\par
\item{40.}
E. Formanek,
The Polynomial Identities and Invariants of $n \times n$ Matrices,
CBMS Regional Conf. Series in Math. 78,
Published for the Confer. Board of the Math. Sci. Washington DC,
AMS, Providence RI, 1991.
\par
\item{41.}
E. Formanek, P. Halpin, W.-C.W. Li,
The Poincar\'e series of the ring of $2 \times 2$ generic matrices,
J. Algebra 69 (1981), 105 -- 112.
\par
\item{42.}
E. Formanek, A.H. Schofield,
Groups acting on the ring of two $2 \times 2$ generic matrices and
a coproduct decomposition of its trace ring,
Proc. Amer. Math. Soc. 95 (1985), 179 -- 183.
\par
\item{43.}
R.M. Guralnick,
Invariants of finite linear groups acting on relatively free
algebras,
Linear Algebra Appl. 72 (1985), 85 -- 92.
\par
\item{44.}
G.B. Gurevich,
Foundations of the Theory of Algebraic Invariants,
P. Noordhoff Ltd., Groningen, 1964.
\par
\item{45.}
H. Helling,
Eine Kennzeichnung von Charakteren auf Gruppen und Assoziativen
Algebren,
Comm. Algebra 1 (1974), 491 -- 501.
\par
\item{46.}
G. Higman,
On a conjecture of Nagata,
Proc. Camb. Philos. Soc. 52 (1956), 1 -- 4.
\par
\item{47.}
D. Hilbert,
\"Uber die Theorie der algebraischen Formen,
Math. Ann. 36 (1890), 473 -- 534;
reprinted in ``Gesammelte Abhandlungen, Band II, Algebra,
Invariantentheorie, Geometrie'', Zweite Auflage,
Springer-Verlag, Berlin - Heidelberg - New York, 1970, 199 -- 257.
\par
\item{48.}
D. Hilbert,
\"Uber die vollen Invariantensysteme,
Math. Ann. 42 (1893), 313 -- 373;
reprinted in ``Gesammelte Abhandlungen, Band II, Algebra,
Invariantentheorie, Geometrie'', Zweite Auflage,
Springer-Verlag, Berlin - Heidelberg - New York, 1970, 287 -- 344.
\par
\item{49.} D. Hilbert,
Mathematische Probleme,
Archiv f. Math. u. Phys. 1 (1901), 44 -- 63, 213 -- 237;
reprinted in ``Gesammelte Abhandlungen, Band III, Analysis,
Grundlagen der Mathematik, Physik, Verschiedenes,
Lebensgeschichte'', Zweite Auflage,
Springer-Verlag, Berlin - Heidelberg - New York, 1970, 290 -- 329.
\par
\item{50.}
N. Jacobson,
PI-Algebras: An Introduction,
Lecture Notes in Math. 441, Springer Verlag, Berlin - New York,
1975.
\par
\item{51.}
I. Kaplansky,
Problems in the theory of rings revised,
Amer. Math. Monthly 77 (1970), 445 -- 454.
\par
\item{52.}
A.R. Kemer,
Capelli identities and nilpotence of the radical of a finitely
generated P.I.-algebra (Russian),
Dokl. Akad. Nauk SSSR 255 (1980), 793 -- 797.
Translation: Soviet Math. Dokl. 22 (1980), 750 -- 753.
\par
\item{53.}
A.R. Kemer,
Asymptotic basis of the identities with unit of the variety ${\rm Var}(M_2(F))$ (Russian),
Izv. Vyssh. Uchebn. Zaved. Mat. (1989), No.6, 71 -- 76.
Translation: Sov. Math. 33 (1990), No.6, 71 -- 76.
\par
\item{54.}
A.R. Kemer,
Ideals of Identities of Associative Algebras,
Translations of Math. Monographs 87, AMS, Providence, RI, 1991.
\par
\item{55.}
V.K. Kharchenko,
Algebra of invariants of free algebras (Russian),
Algebra i Logika 17 (1978), 478 -- 487.
Translation: Algebra and Logic 17 (1978), 316 -- 321.
\par
\item{56.}
V.K. Kharchenko,
Noncommutative invariants of finite groups and Noetherian
varieties,
J. Pure Appl. Algebra 31 (1984), 83 -- 90.
\par
\item{57.}
A.I. Kostrikin,
The Burnside problem (Russian),
Izv. Akad. Nauk SSSR, Ser. Mat. 23 (1959), 3 -- 34.
Translation: Amer. Math. Soc. Transl. (2) 36 (1964), 63 -- 99.
\par
\item{58.}
D. Krakowski, A. Regev,
The polynomial identities of the Grassmann algebra,
Trans. Amer. Math. Soc. 181 (1973), 429 -- 438.
\par
\item{59.}
E.N. Kuz'min,
On the Nagata - Higman theorem (Russian),
in ``Mathematical Structures, Computational Mathematics,
Mathematical Modelling. Proc. Deducated to the 60th Birthday of
Acad. L. Iliev'', Sofia, 1975, 101 -- 107.
\par
\item{60.}
D.R. Lane,
Free Algebras of Rank Two and Their Automorphisms,
Ph.D. Thesis, Bedford College, London, 1976.
\par
\item{61.}
V.N. Latyshev,
Generalization of Hilbert's theorem on the finiteness of bases (Russian),
Sibirsk. Mat. Zh. 7 (1966), 1422 -- 1424.
Translation: Sib. Math. J. 7 (1966), 1112 -- 1113.
\par
\item{62.}
L. Le Bruyn,
The functional equation for Poincar\'e series of trace rings of
generic $2 \times 2$ matrices,
Israel J. Math. 52 (1985), 355 -- 360.
\par
\item{63.}
L. Le Bruyn,
Trace rings of generic 2 by 2 matrices,
Memoirs of AMS, 66, No.363, Providence, R.I., 1987.
\par
\item{64.}
L. Le Bruyn, G. Molenberghs,
Centers of generic division algebras,
Israel Math. Conf. Proc. 1, 1989, 310 -- 319.
\par
\item{65.}
I.V. L'vov,
Maximality cconditions in algebras with identity (Russian),
Algebra i Logika 8 (1969), 449 -- 459.
Translation: Algebra and Logic 8 (1969), 258 -- 263.
\par
\item{66.}
A.I. Malcev,
On the representations of infinite algebras (Russian),
Mat. Sb. 13 (1943), 263 -- 286.
\par
\item{67.}
V.T. Markov,
Systems of generators of T-ideals of finitely generated free algebras (Russian),
Algebra i Logika 18 (1979), 587 -- 598.
Translation: Algebra and Logic 18 (1979), 371 -- 378.
\par
\item{68.}
T. Molien, \"Uber die Invarianten der linearen
Substitutionsgruppen,
Sitz. K\"onig Preuss. Akad. Wiss. (1897), N 52, 1152 -- 1156.
\par
\item{69.} M. Nagata,
On the nilpotency of nil algebras,
J. Math. Soc. Japan 4 (1953), 296 -- 301.
\par
\item{70.}
M. Nagata,
On the 14th problem of Hilbert,
Amer. J. Math. 81 (1959), 766--772.
\par
\item{71.}
E. Noether,
Der Endlichkeitssatz der Invarianten endlicher Gruppen,
Math. Ann. 77 (1916), 89 -- 92;
reprinted in ``Gesammelte Abhandlungen. Collected Papers'',
Springer-Verlag, Berlin - Heidelberg - New York - Tokyo, 1983,
181 -- 184.
\par
\item{72.}
C. Procesi,
Non-commutative affine rings,
Atti Accad. Naz. Lincei, Ser. 8, 8 (1967), 237 -- 255.
\par
\item{73.}
C. Procesi,
Rings with Polynomial Identities,
Marcel Dekker, New York, 1973.
\par
\item{74.}
C. Procesi,
The invariant theory of $n \times n$ matrices,
Adv. in Math. 19 (1976), 306 -- 381.
\par
\item{75.}
C. Procesi,
Trace indentities and standard diagrams,
in ``Ring Theory'', Lect. Notes in Math. 51, Marcel Dekker, New
York, 1979, 191 -- 218.
\par
\item{76.}
C. Procesi,
Computing with $2 \times 2$ matrices,
J. Algebra, 87 (1984), 342 -- 359.
\par
\item{77.}
Yu.P. Razmyslov,
Finite basing of the identities of a matrix algebra of second order
over a field of characteristic 0 (Russian),
Algebra i Logika 12 (1973), 83 -- 113.
Translation: Algebra and Logic 12 (1973), 43 -- 63.
\par
\item{78.}
Yu.P. Razmyslov,
On a problem of Kaplansky (Russian),
Izv. Akad. Nauk SSSR, Ser. Mat. 37 (1973), 483 -- 501.
Translation: Math. USSR Izv. 7 (1973), 479 -- 496.
\par
\item{79.}
Yu.P. Razmyslov,
The Jacobson radical in PI algebras (Russian),
Algebra i Logika 13 (1974), 337 -- 360.
Translation: Algebra and Logic 13 (1974), 192 -- 204.
\par
\item{80.}
Yu.P. Razmyslov,
Trace identities of full matrix algebras over a field of
characteristic zero (Russian),
Izv. Akad. Nauk SSSR, Ser. Mat. 38 (1974), 723 -- 756.
Translation: Math. USSR Izv. 8 (1974), 723 -- 760.
\par
\item{81.}
A. Regev,
Existence of identities in $A \otimes B$,
Israel J. Math. 11 (1972), 131 -- 152.
\par
\item{82.}
L.H. Rowen,
Polynomial Identities in Ring Theory,
Academic Press, New York, 1980.
\par
\item{83.}
G.F. Seelinger,
Generalized matrix valued invariants,
J. Algebra 161 (1993), 199 -- 215.
\par
\item{84.}
G.C. Shephard, J.A. Todd,
Finite unitary reflection groups,
Canad. J. Math. 6 (1954), 274 -- 304.
\par
\item{85.}
K.S. Siberskii,
Algebraic invariants for a set of matrices (Russian),
Sib. Mat. Zh. 9 (1) (1968), 152 -- 164.
Translation: Siber. Math. J. 9 (1968), 115 -- 124.
\par
\item{86.}
W. Specht,
Gesetze in Ringen. I,
Math. Z. 52 (1950), 557 -- 589.
\par
\item{87.}
T.A. Springer,
Invariant Theory,
Lect. Notes in Math. 585, Springer-Verlag, Berlin - Heidelberg -
New York, 1977.
\par
\item{88.}
R.P. Stanley,
Combinatorics and Commutative Algebra,
Progress in Math. 41, Birk\-h\"auser Boston, Boston, MA, 1983.
\par
\item{89.}
Y. Teranishi,
The ring of invariants of matrices,
Nagoya Math. J. 104 (1986), 149 -- 161.
\par
\item{90.}
Y. Teranishi,
Linear diophantine equations and invariant theory of matrices,
``Commut. Algebra and Combinatorics (Kyoto, 1985)'',
Adv. Stud. Pure Math. 11, North-Holland, Amsterdam - New York,
1987, 259 -- 275.
\par
\item{91.}
Y. Teranishi,
The Hilbert series of matrix concominants,
Nagoya Math. J. 111 (1988), 143 -- 156.
\par
\item{92.}
I.K. Tonov,
Just-non-weakly noetherian varieties of associative algebras with unity (Russian),
Pliska Stud. Math. Bulg. 2, 1981, 162 -- 166.
\par
\item{93.}
M. Van den Bergh,
Trace rings are Cohen-Macaulay,
J. Amer. Math. Soc. 2 (1989), 775 -- 799.
\par
\item{94.} M. Van den Bergh,
The center of the generic division algebra,
J. Algebra 127 (1989), 106 -- 126.
\par
\item{95.}
M. Van den Bergh,
Explicit rational forms for the Poincar\'e series of the trace
rings of generic matrices,
Isr. J. Math. 73 (1991), 17 -- 31.
\par
\item{96.}
H. Weyl,
Zur Darstellungstheorie und Invariantenabz\"ahlung der projektiven,
der Kom\-p\-lex- und der Drehungsgruppe,
Acta Math. 48 (1926), 255 -- 278;
reprinted in ``Gesammelte Abhandlungen'', Band III,
Springer-Verlag, Berlin - Heidelberg - New York, 1968, 1 -- 25.
\par
\item{97.}
E.I. Zel'manov,
On Engel Lie algebras (Russian),
Sibirsk. Mat. Zh. 29 (1988), No.5, 112 -- 117.
Translation: Siberian Math. J. 29 (1988), 777 -- 781.

\medskip

\noindent
Vesselin Stoyanov Drensky

\noindent
Institute of Mathematics, Bulgarian Academy of Sciences

\noindent
Akad. G.Bonchev Str., Block 8, 1113 Sofia, Bulgaria

\noindent
e-mail:
drensky@bmath.bas.bg

\bye